\pdfoutput=1
\documentclass[11pt, a4paper, leqno]{article}
\usepackage[intlimits]{amsmath}
\usepackage{amsthm}
\usepackage{amssymb}
\usepackage{graphicx}
\usepackage{fancyhdr}
\numberwithin{equation}{section}
\theoremstyle{plain}% default
\newtheorem{thm}{Theorem}[section]
\newtheorem*{algorithm}{Main Algorithm}

\newtheorem{lem}[thm]{Lemma}

\theoremstyle{definition}

\theoremstyle{remark}
\newtheorem{rem}[thm]{Remark}

\title{A simple method of calculating eigenvalues and resonances in
domains with infinite regular ends\thanks{To appear in {\em Proc. Royal Soc. Edinb.} (2008)}}
\DeclareMathOperator{\im}{Im}
\DeclareMathOperator{\re}{Re}

\DeclareMathOperator{\supp}{supp}
\DeclareMathOperator{\spec}{spec}

\DeclareMathOperator{\diver}{div}
\newcommand{\bgrad}{\text{\bf grad}\,}
\renewcommand{\Im}{\im}
\renewcommand{\Re}{\re}
\newcommand{\yb}{\mathbf{y}}
\newcommand{\xb}{\mathbf{x}}

\newcommand{\Rcal}{\mathcal{R}}
\newcommand{\Ccal}{\mathcal{C}}
\newcommand{\Acal}{\mathcal{A}}

\newcommand{\Pcal}{\mathcal{P}}
\newcommand{\Qcal}{\mathcal{Q}}

\newcommand{\Stt}{\mathsf{S}}
\newcommand{\Mtt}{\mathsf{M}}
\newcommand{\Dtt}{\mathsf{D}}
\newcommand{\Ttt}{\mathsf{T}}
\newcommand{\Rtt}{\mathsf{R}}
\newcommand{\cvect}{\mathsf{c}}
\newcommand{\zerovect}{\mathsf{0}}
\newcommand{\Rbb}[1][]{\mathbb{R}^{#1}}

\newcommand{\Nbb}{\mathbb{N}}
\newcommand{\dpar}{\partial}

\newcommand{\partd}[2]{{\frac{\dpar #1}{\dpar #2}}}
\newcommand{\scal}[2]{{\langle#1,#2\rangle}}
\author{Michael Levitin  and Marco Marletta%
\\
\normalsize\small School of Mathematics,  Cardiff University,\\
\normalsize\small Senghennydd Road, Cardiff CF24 4AG, UK\\
\normalsize\small email {\sffamily \{Levitin, MarlettaM\}@cardiff.ac.uk}\\
\normalsize\small {\sffamily http://www.cf.ac.uk/maths/subsites/levitin, http://www.cf.ac.uk/maths/subsites/marlettam/ }
}
\begin{document}
%
%
%\subjclass{}
%
\date{\small Final version --- 20th August 2008}
%\\\fbox{\bf Preliminary version - not for distribution}}
%
\maketitle
\begin{abstract} \noindent A new approach is presented for the solution of spectral 
problems on infinite domains with regular ends, which avoids the need to solve boundary value problems for many trial values of the spectral parameter. We present numerical results both for eigenvalues and for resonances, comparing with results reported by Aslanyan, Parnovski and Vassiliev [{\em Q. J. Mech. Appl. Math.} {\bfseries 53} (2000), 429--447].
\end{abstract}

\

%%%%%%%%%%%%%%%%%%%%%%%%%%%%%%%%%%%%%%%%
\newpage
\tableofcontents
\newpage
%%%%%%%%%%%%%%%%%%%%%%%%%%%%%%%%%%%%%%%%%
\section{Introduction}
We present a simple new approach to the solution of a wide class of spectral and resonance
problems on infinite domains with regular ends, including those found in the study of 
quantum switches, waveguides, and acoustic scatterers. Our algorithm is part analytical
and part numerical and is essentially a combination of four classical approaches
(domain decomposition, boundary elements, finite elements and spectral methods) each
of which is used in its most natural context. 

Our method has three main advantages. Firstly, it allows one to exploit existing library 
software for PDEs to reduce to a minimum the amount of new code which must 
be written: meshing, choice of basis functions, preconditioning and solution 
of linear systems can all be delegated to the library code. Secondly, the method 
is well adapted to handling resonance or scattering problems, where Neumann-to-Dirichlet 
maps must be continued analytically as functions of the spectral 
parameter. Thirdly, the need to solve boundary value problems for many different 
choices of the spectral parameter is totally removed.

A brief outline of the solution procedure is as follows. Firstly, the infinite domain
on which the problem is posed is decomposed ({\em domain decomposition}) into a union
of a bounded domain and a finite number of infinite regular ends, e.g. cylinders. We
need to assume that on each of these ends, the partial differential equation admits
separation of variables. The problem can thereby be reduced to a problem on the
finite portion of the domain with extra boundary conditions on the interfaces with
the attached ends. This problem can be reduced to a problem on the interfaces through
the interior `Neumann-to-Dirichlet' map, whose matrix elements, with respect to a
specially chosen basis on the interfaces ({\em boundary element method}) are calculated 
in terms of traces of `mixed' Neumann eigenfunctions from the interior domain ({\em spectral
method}), which is probably the main new idea of this article.
The Neumann eigenfunctions are themselves computed by a {\em finite element method}.
The resulting problem on the interface can be solved by finding roots of some monotone
functions calculated by solving a generalized matrix eigenproblem of modest
size. This is done by elementary means.

Apart from the need to be able to separate variables at infinity, the remaining parts
of our approach work for a wide class of problems with minimal restrictions on the
PDE coefficients and the boundary conditions. For simplicity of exposition, however,
most of this paper is presented for the case of the Helmholtz equation in domains
with cylindrical ends, or in domains exterior to an obstacle.

Since the 1960s there has been a lot of important theoretical and numerical
work done on the solution of spectral and scattering problems on the sort of
domains we consider. The literature is too extensive to review adequately
here: we mention just a few trends.

On the purely mathematical side, attention has recently concentrated on finding 
conditions sufficient for existence of eigenvalues either below the 
essential spectrum or embedded in it, see e.g. \cite{EvLeVa, ES, DE, KK}. 
Additionally, there are some partial results on non-existence of eigenvalues 
in the lower part of the essential spectrum \cite{DaPa}, and on perturbations 
of embedded eigenvalues which turn  them into resonances. For an abstract approach 
to counting the eigenvalues, see \cite{Pa, Ch, ChZw}.

For numerical treatments, boundary integral and boundary element methods have
often been employed --- see, e.g., Colton and Kress \cite{Colton}, Brebbia,
Telles and Wrobel \cite{Brebbia} for reviews. Others have considered
domain truncation methods coupled with the use of `non-reflecting' 
boundary conditions --- see, e.g., H\"{a}gstrom and Keller \cite{Hagstrom};
Givoli, Patlashenko and Keller \cite{Givoli}. For the particular case
of waveguides, special efficient numerical approaches can be found in Aslanyan, 
Parnovski and Vassiliev \cite{APV} and in McIver, Linton, McIver and 
Zhang \cite{McIver}; the scattering matrix approach has also been developed
in, e.g., Bagraev, Mikhailova, Pavlov, Prokhorov and Yafyasov \cite{BMPPY} and
in Semenikhin, Pavlov and Ryzhii \cite{SPR}. Various related numerical methods 
and applications to a variety of physically motivated special cases and generalisations 
can be found, e.g., in \cite{Ur, BoJo, HaPaGi, PeJo, Bo} and references therein. 

We do not claim that our approach is more numerically efficient than any of these: 
what is certainly true, however, is that it is much simpler to implement, requiring 
just a few lines of MATLAB code rather than some serious work by expert programmers.

Some examples of our code are available from the authors' web pages; we encourage the reader
to download, modify, and use it.

\section{Geometry}\label{sec:geom}
By $\Omega\subset\Rbb[d]$, $d\ge 2$ we denote an unbounded open connected domain with $N$ non-intersecting
cylindrical ends. Namely, for each $1\leq n \leq N$ we choose 
local coordinates $(x,\yb)$ such that the cylinder $\Ccal_n^0$ is given by
\[ 
\Ccal_n^0 = \{ (x,\yb) \; : \; \yb\in \Gamma_n^0, \;\; x \ge 0\}, 
\]
where $\Gamma_n^0$ is the open bounded connected (but not necessarily simply connected) 
$(d-1)$-dimensional cross-section of $\Ccal_n^0$, and we suppose that
\[
 \overline{\Omega}=\overline{\Omega_0\sqcup\bigsqcup_{n=1}^N \Ccal_n^0}\, ,
\]
where $\Omega_0\subset\Rbb[d]$ is a connected bounded open set. Further geometrical definitions
(interior part of $\Omega$, cylindrical ends of $\Omega$ and interfaces) are introduced in the 
next section.

Throughout the paper we assume that the boundary $\partial\Omega$ is at least piecewise $C^1$ smooth and, for simplicity,
satisfies both interior and exterior uniform cone conditions (i.e.  $\partial\Omega$ does not have any cusps).

\section{Spectral problem; boundary conditions; additional notation}\label{sec:sp_prob}

We consider, in $\Omega$, the spectral problem
%%%%
\begin{equation}\label{eq:lapl}
-\Delta u = \lambda u
\end{equation}
%%%%
subject to the boundary conditions 
%%%%
\begin{equation}\label{eq:bc}
B u:= \left.\left(a(\cdot)u+b(\cdot)\partd{u}{n}\right)\right|_{\dpar\Omega} = 0\,.
\end{equation}
Here the unit normal $n$ to $\dpar\Omega$ points outwards.
%%%%

We shall refer to Problem \eqref{eq:lapl}, \eqref{eq:bc} as ``Problem (P)''.

We impose the following restrictions on the boundary coefficients $a(\cdot)$ and $b(\cdot)$:
\begin{itemize}
\item[(B$_1$)] $a$ and $b$ are piecewise smooth functions on the boundary $\dpar\Omega$;
\item[(B$_2$)] $a(\cdot)^2+b(\cdot)^2\equiv 1$;
\item[(B$_3$)] on each connected component of the part $(\frac{1}{2},\infty)\times \partial\Gamma_n$ 
of the infinite boundary  of each cylindrical end $\Ccal_n^0$ the functions $a$ and $b$ are constant 
(but the constants are allowed to differ between the cylinders and even between the connected components 
of the boundary of each cylinder).
\end{itemize}

To finalize our notation, we shall from now on refer to the sets
\[
\Ccal_n:=\{ (x,\yb) \; : \; \yb\in \Gamma_n^0, \;\; x \ge 1\}
\]
and their union 
\[
\Ccal:=\bigsqcup_{n=1}^N \Ccal_n
\]
as \emph{cylindrical ends} of $\Omega$ and the sets 
\[
\Gamma_n:=\{(1,\yb)\in\Ccal_n\}=\{ (1,\yb) \; : \; \yb\in \Gamma_n^0\}
\]
and their union
\[
\Gamma:=\bigsqcup_{n=1}^N \Gamma_n
\]
as \emph{interfaces}.

\begin{figure}[htb!]
\begin{center}
\resizebox{260pt}{170pt}{\includegraphics*{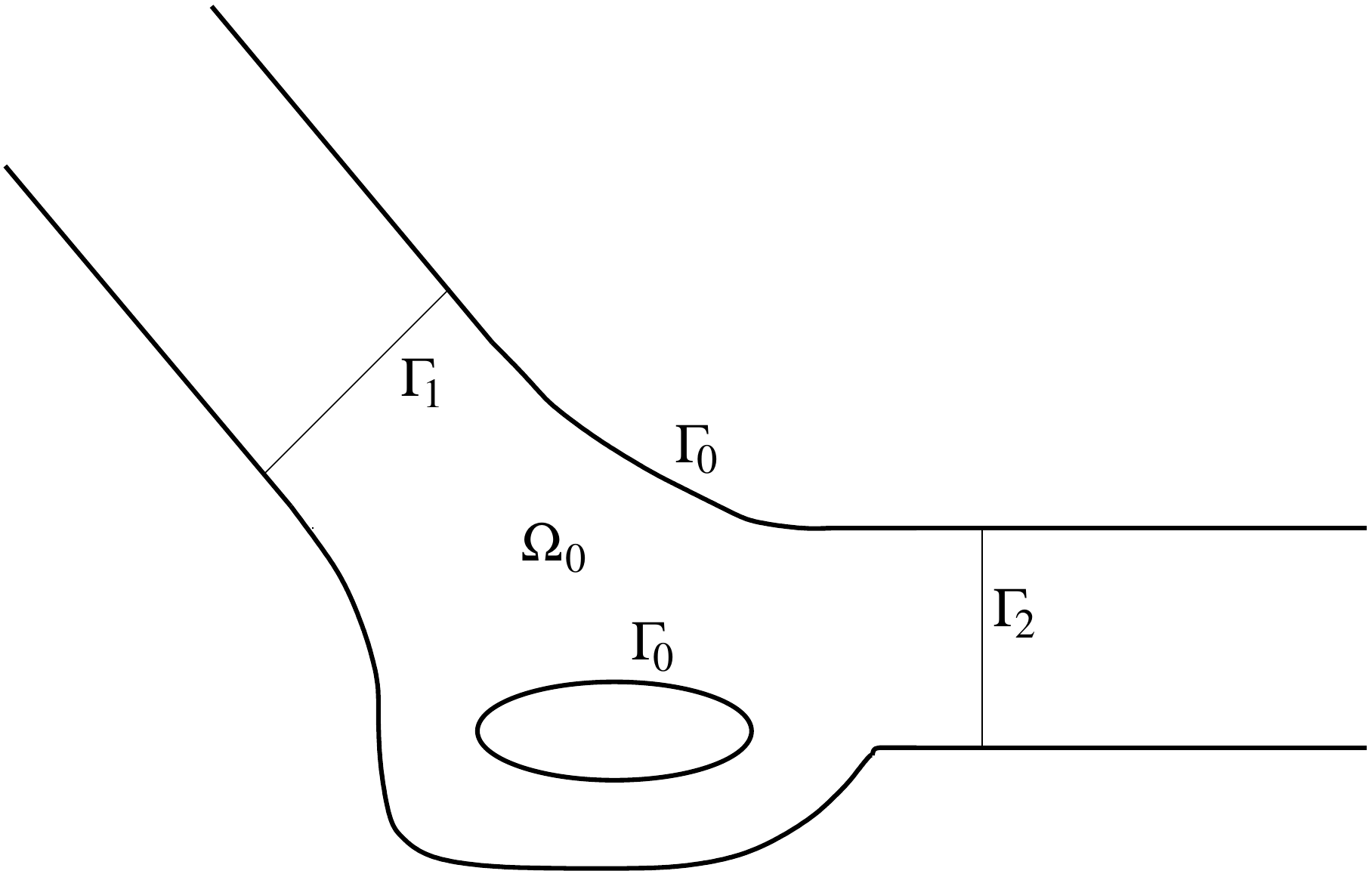}}
\end{center}
\caption{Example of a domain with cylindrical ends.\label{fig:domain}}
\end{figure}

We also define the \emph{interior part} of $\Omega$ as
\[
\Omega_0 := \Omega\setminus \Ccal
\]
and denote by 
\[
\Gamma_0 := \dpar\Omega_0\setminus\Gamma\,,
\]
the part of its boundary shared with the boundary of $\Omega$, see Figure \ref{fig:domain}.

%%%%%%%%%%%%%%%%%%%%%
%%%%%%%%%%%%%%%%%%%%%%%

We emphasize once more  that according to (B$_3$) the boundary condition \eqref{eq:bc} are constant on 
each connected component of 
\[
\dpar\Ccal\setminus\Gamma = \bigsqcup_{n=1}^N (1,\infty)\times\dpar\Gamma_n\,.
\]

\section{Generalized problem}\label{sec:genprob}

As the reader shall easily see later on, all our methods apply, with minimal modifications, to the general second order spectral boundary value
problem
\begin{equation}\label{eq:gen_sp_prob}
-\diver(p(\cdot)\bgrad u)+q(\cdot) u = \lambda u\,
\end{equation} 
\begin{equation}\label{eq:gen_bc}
\left.\left(a(\cdot)u+b(\cdot)\mathbf{n}\cdot p(\cdot)\bgrad u(\cdot)\right)\right|_{\dpar\Omega} = 0\,.
\end{equation}
with $p$ being a sufficiently smooth positive definite $(d\times d)$-matrix valued function, and $q$ being an $L_\infty(\Rbb[d])$ scalar potential.
Then, in addition to assumptions (B${}_1$)--(B${}_3$), we also require that $p(\cdot)$ and $q(\cdot)$ be constant on each cylindrical end 
$\Ccal_n^0$  (but the constants are allowed to differ between the cylinders).

Although one can easily formulate all the results and algorithms for the general case \eqref{eq:gen_sp_prob}, \eqref{eq:gen_bc}, we deliberately describe only 
Problem (P) for the Laplacian in order to preserve the simplicity and clarity of exposition.

\section{Rigorous problem statement}\label{sec:rig} 
We always understand the spectral problem (P) in the variational sense. Namely, let 
$\partial_D\Omega := \{\xb\in\partial\Omega\;:\; b(\xb)=0\}$ be the part of the boundary on 
which \eqref{eq:bc} is the Dirichlet condition. 
Let $C_0^\infty(\Omega):=\{v\in C^\infty(\Omega)\;:\;\overline{\supp v}\cap\partial_D\Omega=\emptyset\}$ be 
the set of all infinitely differentiable functions on $\Omega$ supported away from $\partial_D\Omega$. 
We define the sesquilinear form 
\begin{equation}\label{eq:form}
\mathfrak{t}[u,v]:=\int_\Omega \nabla u\cdot\overline{\nabla v} + \int_{\partial\Omega\setminus \partial_D\Omega} 
\frac{a}{b} u\overline{v}
\end{equation}
first on $C_0^\infty(\Omega)$, and then extend it to the closure $H^1_*(\Omega)$ of $C_0^\infty(\Omega)$ with respect to the 
norm associated with \eqref{eq:form}.

From now on, we shall understand Problem (P) as the spectral problem $\mathfrak{t}[u,v]=\lambda\scal{u}{v}$, where 
$\scal{\cdot}{\cdot}$ denotes the scalar product in $L_2(\Omega)$, with the form domain in $H^1_*(\Omega)$. In particular 
$\lambda$ is an \emph{eigenvalue} of Problem (P) if there exists a nontrivial function $u\in H^1_*(\Omega)$ such that 
${\mathfrak t}(u,v) = \lambda \scal{u}{v}$ for all $v\in H^1_*(\Omega)$.

\section{Transversal problem and essential spectrum}\label{sec:trans} 

We need to consider the following \emph{transversal} problem on the joint interface $\Gamma$:
\begin{equation}\label{eq:trans}
-\Delta_{d-1} w = \kappa w\,,\qquad \left.\left(a(\cdot)w + b(\cdot)\partd{w}{\nu}\right)\right|_{\partial\Gamma} =0\,;
\end{equation}
here $\Delta_{d-1}$ is the $(d-1)$-dimensional Laplacian and $\partd{}{\nu}$ is the derivative with respect to the external 
normal $\nu$ to  $\partial\Gamma$ (which of course coincides here with the normal $n$ to $\partial\Omega$). We recall
that $a$ and $b$ are constant on each connected component of $\partial\Gamma$. 

This problem has a rigorous formulation in terms of sesquilinear forms analogous to that introduced in (\ref{eq:form})
but with $\Omega$ replaced by $\Gamma$ throughout. The assumptions (B$_1$), (B$_2$) and (B$_3$)
ensure that the problem \eqref{eq:trans} has a purely discrete spectrum of eigenvalues $(\kappa_j)_{j=1}^\infty$, enumerated in increasing order and repeated according to multiplicity; we denote by $(w_j(\yb) )_{j=1}^\infty$ the corresponding eigenfunctions normalized in $L_2(\Gamma)$. We also set $\kappa_0:=-\infty$.

It is well known that the essential spectrum of Problem (P) is then given by $[\kappa_1,+\infty)$ and, moreover, that its
multiplicity at a point $\lambda\ge\kappa_1$ is equal to the number of $\kappa_j$ which are less than or equal to $\lambda$. We shall often refer to the number $\kappa_j$ as the $j$-th \emph{threshold}. 

For each $J\geq 1$ we define
\[
P_J := \mbox{Span}\{w_j\, : \, j > J \}\,,
\]
set $P_0 := L_2(\Gamma)$, and denote by $\Pcal_J$ the projection in $L_2(\Gamma)$ onto $P_J$ ($J\ge 0$). We also define $\Qcal_J := I-\Pcal_J$.

\section{Eigenvalue problem --- reduction to the interface}\label{sec:reduct}
We define formally two Neumann to Dirichlet operators, for the interior part $\Omega_0$ and
cylindrical ends $\Ccal$. 

Firstly, the interior Neumann to Dirichlet operator $\Rcal_{\lambda}^0:=g\mapsto \left.v\right|_\Gamma$
is the map which maps any suitable function $g$ defined on $\Gamma$ to the trace on
$\Gamma$ of the solution of the boundary value problem
\begin{equation}\label{eq:intg}
-\Delta v = \lambda v \;\; \text{in} \;\; \Omega_0, \quad 
\frac{\partial v}{\partial n}=g  \;\; \text{on} \;\;\Gamma, \quad Bv = 0 \;\; \text{on} \;\; \Gamma_0 
\end{equation}
(where by $n$ we always mean the unit normal to $\Gamma$ directed outwards from $\Omega_0$)
provided $\lambda$ is not an eigenvalue for the corresponding Neumann problem
\begin{equation}\label{eq:intN}
 -\Delta U = \lambda U \;\; \text{in} \;\; \Omega_0, \quad \frac{\partial U}{\partial n}=0  \;\; \text{on} \;\;
 \Gamma, \quad BU = 0 \;\; \text{on} \;\; \Gamma_0.
\end{equation}
 
Secondly, the exterior Neumann to Dirichlet operator $\Rcal_{\lambda}^{\Ccal}:=g\mapsto \left.V\right|_\Gamma$
is the map which maps any suitable function $g$ defined on $\Gamma$ to the trace on
$\Gamma$ of the solution of the exterior boundary value problem
\begin{equation}
-\Delta V = \lambda V \;\; \text{in} \;\; \Ccal, \quad -\frac{\partial V}{\partial n}=g  \;\; \text{on} \;\;
 \Gamma, \quad BV = 0 \;\; \text{on} \;\; \partial\Ccal\setminus\Gamma 
\label{eq:extNeum}
\end{equation}
(the normal $n$ on $\Gamma$, as before, points outwards from $\Omega_0$ along the positive $x$-axis of each cylindrical end, and therefore towards $\Ccal$).

Suppose that $u$ is an eigenfunction of Problem (P) with eigenvalue $\lambda$. It is well known
that $u$ is infinitely differentiable at any interior point of $\Omega$. In particular, therefore, 
$u$ and its gradient restrict to the interface $\Gamma$ as smooth functions, and so we may
define
\begin{equation}\label{eq:fg}
f := \left. u \right|_{\Gamma}, \quad g = \left. \frac{\partial u}{\partial n}\right|_{\Gamma}. 
\end{equation}
Thus $f$ and $g$ satisfy 
\[ 
f = -\Rcal_\lambda^{\Ccal}g = \Rcal_\lambda^0g, 
\]
so that a necessary condition for $\lambda$ to be an eigenvalue of Problem (P) is that $\sigma = -1$
be an eigenvalue
of the operator pencil problem
\begin{equation}\label{eq:pencil} 
(\sigma\Rcal_\lambda^\Ccal-\Rcal_\lambda^0)g = 0\, .
\end{equation}
(A number $\sigma_0$ is called an eigenvalue of a general operator pencil $\sigma A + B$ if  zero is an 
eigenvalue of the operator $\sigma_0 A + B$.)

%A graph of a typical behaviour of one of the eigenvalues $\sigma$ of (\ref{eq:pencil})
%as a function of $\lambda$ is shown in Fig. \ref{sig}. In fact, this is one of the eigenvalues
%of the pencil which arises for a bent waveguide problem --- Example 1 in Section 
%\ref{numalg}.
%%%%%%%%%%%%%%%%%%%%%
%\begin{figure}
%\begin{center}
%\resizebox{250pt}{220pt}{\includegraphics*{sigma1}}
%\end{center}
%\caption{Pencil eigenvalue as a function of $\lambda$\label{sig}}
%\end{figure}
%%%%%%%%%%%%%%%%%%%%%%%

\section{Problem on a cylindrical end --- the Neumann to Dirichlet operator}\label{sec:NDcyl}
Suppose that $V$ is a solution of the exterior problem (\ref{eq:extNeum}) for some $\lambda$ (not necessarily
an eigenvalue of Problem (P)). A formal expansion of $V$ in terms of the eigenfunctions $w_j$ of the transversal
problem (\ref{eq:trans}), 
\[ 
V(x,\yb) = \sum_{j=1}^\infty c_j(x) w_j(\yb) 
\]
with account of \eqref{eq:extNeum} yields for the functions $c_j$ the ordinary differential boundary value problems 
\begin{equation}\label{eq:cj}
\frac{d^2 c_j(x)}{dx^2}+(\lambda-\kappa_j)c_j(x)=0\,\quad\text{for }x>1\,,\qquad 
\left.\frac{dc_j(x)}{dx}\right|_{x=1}=-(g,w_j)\,.
\end{equation}
where $(\cdot,\cdot)$ denotes the scalar product in $L_2(\Gamma)$.

The two solutions of \eqref{eq:cj} are
\begin{equation}\label{eq:cjpm}
c_j^\pm(x)=\mp\frac{(g,w_j)}{\sqrt{\kappa_j-\lambda}}\exp\left(\pm (x-1)\sqrt{\kappa_j-\lambda}\right)\,.
\end{equation}

If we additionally require that the solution $V\in H^1_*(\Ccal)$ and assume that, for some integer $J\ge0$, 
$\lambda\in[\kappa_J,\kappa_{J+1})$, then the boundary Neumann data $g$ should belong to the space $\Qcal_J H^{-1/2}_*(\Gamma)$ 
(in particular, we need $(g,w_j)=0$ for $j\le J$)  and
choose only the solutions $c_j(x):=c_j^-(x)$ in order to to
exclude terms growing or oscillating at $x=+\infty$. Then
$V(1, \yb) = \sum_{j=J+1}^\infty  c_j(1)w_j(\yb)\in H^{1/2}_*(\Gamma)$ and we have therefore defined the non-negative self-adjoint operators
\[
\Rcal_\lambda^\Ccal: \Qcal_J H^{-1/2}_*(\Gamma)\to \Qcal_J H^{1/2}_*(\Gamma)\quad\text{for }\lambda\in[\kappa_J,\kappa_{J+1})
\]
acting as
\begin{equation}\label{eq:rlamC}
\Rcal_\lambda^\Ccal g = \sum_{j=J+1}^\infty \frac{1}{\sqrt{\kappa_j-\lambda}}(g,w_j)w_j.
\end{equation}

It is clear from \eqref{eq:rlamC} and from the construction  
that the quadratic form $(\Rcal_\lambda^\Ccal g, g)_\Gamma$ is monotone increasing in $\lambda$ on each interval $[\kappa_J,\kappa_{J+1})$.

\section{Problem on an interior domain --- the Neumann to Dirichlet operator}\label{sec:NDint}
Let $g_1$ and $g_2$ be functions defined on $\Gamma$ and suppose that $v_1$ and $v_2$ are the
corresponding solutions of the boundary value problem (\ref{eq:intg}), for $\lambda$ which is not an
eigenvalue of (\ref{eq:intN}). A formal integration by parts shows that
\begin{equation}\label{eq:qform}
(\Rcal_\lambda^0 g_1, g_2) = (g_1,\Rcal_\lambda^0 g_2) = \langle \nabla v_1, \nabla v_2 \rangle
 - \lambda \langle v_1, v_2 \rangle\, . 
\end{equation}

This allows us to define a quadratic form 
\[
\mathfrak{r}_\lambda[g]:=(\Rcal_\lambda^0 g,g) = \langle \nabla v, \nabla v \rangle
 - \lambda \langle v, v \rangle
\]
(where $v$ solves \eqref{eq:intg}) associated with $\Rcal_\lambda^0$
and hence to regard $\Rcal_\lambda^0$ either as a self-adjoint operator 
on a domain in $L_2(\Gamma)$ or, under suitable 
smoothness restrictions on the boundary $\partial\Omega_0$, as a map 
on a scale of Sobolev spaces from $H^{s-1/2}(\Gamma)$ to $H^{s+1/2}(\Gamma)$ 
(see \cite{Agranovich,Paltsev,Safarov}). We emphasize 
that we impose a condition that $\lambda$ is not an eigenvalue of (\ref{eq:intN}) --- 
for ways of bypassing this condition see \cite{Safarov}.

For a fixed $g\in L_2(\Gamma)$, the form $\mathfrak{r}_\lambda[g]$ is known to be monotone 
increasing with respect to $\lambda$ in any interval not containing the eigenvalues
of (\ref{eq:intN}), as the following simple argument (found e.g. in \cite{Friedlander}) 
shows. Let $g$ be fixed and let $v'$ denote the derivative with respect to $\lambda$ of 
the solution $v$ of (\ref{eq:intg}). Then $v'$ solves the boundary value problem
\begin{equation}\label{eq:vprime}
-\Delta v' = \lambda v' + v\quad\text{in }\Omega_0\,,\quad \partial v'/\partial n|_\Gamma=0\,,\quad Bv'|_{\Gamma_0}=0\,.
\end{equation} 
Multiplying \eqref{eq:vprime} by $\overline{v}'$ and integrating by parts we obtain
\begin{equation}\label{eq:monot}
\frac{1}{2}\,\frac{d}{d\lambda}\mathfrak{r}_\lambda[g] = \left(\frac{d}{d\lambda}(\Rcal^0_\lambda g), g\right)_\Gamma = 
\left(v', \frac{\partial v}{\partial\nu}\right)_\Gamma = \langle v, v \rangle_{\Omega_0} \ge 0,
\end{equation} 
which proves the claim. Exactly the same argument can be used to prove the monotonicity in $\lambda$ of the form of the operator 
$\Rcal_\lambda^\Ccal$, which we established at the end of the previous section by direct inspection of \eqref{eq:rlamC}.

\section{Problem on an interior domain --- basis expansions}\label{sec:NDexp}

In \S~\ref{sec:NDcyl}, we constructed an explicit representation \eqref{eq:rlamC} for the Neumann to Dirichlet operator 
$\Rcal_\lambda^\Ccal$ associated with the problem on a cylindrical end. The aim of this section is to do the same for the 
\emph{interior} Neumann to Dirichlet map $\Rcal_\lambda^0$, and we want to emphasize that the following trick is in fact the 
focal point of this paper.

In the setting of problem \eqref{eq:intg}, let us choose an \emph{arbitrary} basis $\left\{\phi_k\right\}_{k=1}^\infty$ in 
$L_2(\Gamma)$ (which is the domain of the quadratic form $\mathfrak{r}_\lambda$).
Let us also denote by $\Phi_k$ the solution $v$ of \eqref{eq:intg} with $g=\phi_k$, i.e. of
\begin{equation}\label{eq:intPhi}
-\Delta \Phi_k = \lambda \Phi_k \;\; \text{in} \;\; \Omega_0, \quad 
\frac{\partial \Phi_k}{\partial n}=\phi_k  \;\; \text{on} \;\;\Gamma, \quad B\Phi_k = 0 \;\; \text{on} \;\; \Gamma_0. 
\end{equation}

We also denote by $\mu_m$ the eigenvalues and by $U_m$ the corresponding eigenfunctions of the homogeneous problem \eqref{eq:intN}, i.e.
\begin{equation}\label{eq:U}
 -\Delta U_m = \mu_m U_m \;\; \text{in} \;\; \Omega_0, \quad \frac{\partial U_m}{\partial n}=0  \;\; \text{on} \;\;
 \Gamma, \quad BU_m = 0 \;\; \text{on} \;\; \Gamma_0, 
\end{equation}

Our aim is to find explicit expressions for the matrix elements $R_{k\ell}:=(\Rcal_\lambda^0 \phi_k, \phi_\ell)_\Gamma$. 

By \eqref{eq:qform} we have
\begin{equation}\label{eq:startexp}
R_{k\ell}=(\Rcal_\lambda^0 \phi_k, \phi_\ell)_\Gamma=\langle\nabla\Phi_k,\nabla\Phi_\ell\rangle-\lambda\langle\Phi_k,\Phi_\ell\rangle\,.
\end{equation}

We now use the fact that the eigenfunctions $U_m$ form a basis in the Hilbert space $H^1_*(\Omega_0)$. We can thus expand $\Phi_k$
in this basis as 
\[
\Phi_k=\sum_{m=1}^\infty U_m\langle\Phi_k, U_m\rangle\,. 
\]
Thus , expanding the right-hand side of \eqref{eq:startexp} with the use of two elementary integrations by parts, we obtain
\begin{equation}\label{eq:mainexpr}
\begin{split}
R_{k\ell}&=\sum_{m=1}^\infty {\langle\nabla U_m, \nabla U_m\rangle}{\langle\Phi_k,U_m\rangle}{\langle U_m,\Phi_\ell\rangle}\\
&=\sum_{m=1}^\infty (\mu_m-\lambda){\langle\Phi_k,U_m\rangle}{\langle U_m,\Phi_\ell\rangle}
=\sum_{m=1}^\infty \frac{1}{\mu_m-\lambda}(\phi_k,U_m|_\Gamma)_\Gamma\cdot(U_m|_\Gamma,\phi_\ell)_\Gamma\,.
\end{split}
\end{equation}

\section{Embedded eigenvalues and orthogonality conditions}
We have already established formally in \S~\ref{sec:reduct} that if $\lambda$ is an eigenvalue of Problem (P) then 
$\sigma=-1$ is an eigenvalue of the  pencil  \eqref{eq:pencil} on the interface $\Gamma$.
In view of rigorous definitions of the exterior Neumann to Dirichlet operator, we 
want now to consider, for $\lambda\in[\kappa_J, \kappa_{J+1})$, the pencil 
\begin{equation}\label{eq:pencildef}
\Acal_\lambda(\sigma) := \sigma\Rcal^\Ccal_\lambda-\Qcal_J\Rcal^0_\lambda\,.
\end{equation} 
acting in a weak sense in $\Qcal_J L_2(\Omega)$.

For each fixed $\lambda$ (which plays the role of a parameter rather than spectral parameter) $\Acal_\lambda(\sigma)$ has $\sigma$-spectrum
$\spec(\Acal_\lambda)$, and a necessary condition for $\lambda$ to be an eigenvalue of Problem (P) is that $-1\in\spec(\Acal_\lambda)$.
The corresponding eigenfunction $g\in\Qcal_JL_2(\Gamma)$ is the normal derivative on $\Gamma$ of the eigenfunction $u$ of 
Problem (P) corresponding to $\lambda$.
Also, if $\lambda\in [\kappa_{J},\kappa_{J+1})$ for some $J>0$ then we have seen in \S~\ref{sec:NDcyl} that 
$g$ must additionally satisfy the orthogonality condition $\Pcal_J\Rcal^0_{\lambda}g = (I-\Qcal_J)\Rcal^0_{\lambda}g=0$; in other words
$\Rcal^0_{\lambda}g\in \Qcal_JL_2(\Omega)$. 

Since $\Qcal_J L_2(\Omega)$ is an invariant space of the operator
$\Rcal_\lambda^\Ccal$, we propose the following strategy for finding embedded eigenvalues:
\begin{itemize}
\item first find those ``suspicious points'' $\lambda_*$ for which $\sigma=-1$ is an eigenvalue of the operator pencil  \eqref{eq:pencildef}
 on the interface $\Gamma$;
\item if these ``suspicious'' $\lambda_*$ lie in $[\kappa_J,\kappa_{J+1})$ for some $J>0$, check the orthogonality
 condition $\Pcal_J\Rcal^0_{\lambda_*} g=0$ for the corresponding pencil eigenfunction $g\in\Qcal_JL_2(\Gamma)$.
\end{itemize}

Assume that we wish to find the eigenvalues in the interval $[\kappa_J,\kappa_{J+1})$. 
We start by reducing the spectral problem for the pencil \eqref{eq:pencildef} on the space $\Qcal_J L_2(\Gamma)$ to a spectral problem for an 
infinite matrix pencil. In \S~\ref{sec:NDcyl} we developed a representation of the operator $\Rcal_\lambda^\Ccal$ in 
terms of the transverse eigenfunctions $w_j$; in \S~\ref{sec:NDexp} we developed a representation for $\Rcal_\lambda^0$ in terms
of an arbitrary basis $\{ \phi_k \}_{k=1}^\infty$ in $L_2(\Gamma)$. We now choose as our basis of $\Qcal_J L_2(\Gamma)$ the
set $\{ \phi_k \}_{k=1}^\infty$ given by $\phi_k = w_{k+J}$, $k\in \Nbb$.
Consequently we obtain from \eqref{eq:mainexpr} that
\begin{equation}\label{eq:Rklinner}
 R_{k\ell} := (\Rcal_\lambda^0w_{k+J},w_{\ell+J})_\Gamma = \sum_{m=1}^\infty \frac{1}{\mu_m-\lambda}(w_{k+J},U_m|_\Gamma)_\Gamma
 \cdot (U_m|_\Gamma,w_{\ell+J})_\Gamma.
\end{equation}
We simplify the expression for the infinite matrix $\Rtt(\lambda)=(R_{k\ell})_{1\leq k,\ell<\infty}$ by introducing some notation. 
Define an infinite matrix $\Stt$ 
with entries $\left(S_{km}\right)_{1\leq k,m <\infty}$ defined by
\begin{equation}\label{eq:Sdef} 
 S_{km} = \left(w_k,U_m|_\Gamma\right)_\Gamma,
\end{equation}
an infinite diagonal matrix $\Dtt(\lambda)$ with diagonal entries 
\begin{equation}\label{eq:Ddef}
 D_{kk}(\lambda) = (\mu_k-\lambda)^{-1}, \;\;\; 1\leq k < \infty,
\end{equation}
and an infinite diagonal matrix $\Ttt(\lambda)$ with diagonal entries
\begin{equation}\label{eq:Tdef}
 T_{kk}(\lambda) = \frac{1}{\sqrt{\kappa_k-\lambda}}, \;\;\; 1\leq k <\infty. 
\end{equation}
For a general matrix $\Mtt$ we will use the notation $M_{a:b,c:d}$ to denote the submatrix 
\begin{equation}\label{eq:Mprop}
    \left(M_{ij}\right)_{{a\leq i \leq b}\atop{c\leq j \leq d}},
\end{equation}
where $b=\infty$ and $d=\infty$ are allowed.

The matrix $\Rtt$ whose entries are defined in \eqref{eq:Rklinner} can now be written as
\begin{equation} \label{eq:matform}
 \Rtt(\lambda) = \Stt_{J+1:\infty,1:\infty} \Dtt(\lambda) \left(\Stt_{J+1:\infty,1:\infty}\right)^*\, .
\end{equation}

In a similar way, we can define the matrix elements associated with $\Rcal_\lambda^\Ccal$ with
respect to the same basis using \eqref{eq:rlamC}: $\Rcal_\lambda^\Ccal$ has a diagonal
matrix representation $\Ttt(\lambda)$ in this basis.
We thus obtain the following reduction of the problem from $\Qcal_J L_2(\Gamma)$ to $\ell^2(\Nbb)$.

\begin{thm}\label{thm:equiv}
A number $\lambda\in[\kappa_J,\kappa_{J+1})$ (which is assumed not to be an eigenvalue of the 
auxilliary Neumann problem \eqref{eq:U}) is an eigenvalue of Problem (P) if and only if $-1$ is a $\sigma$-eigenvalue 
of the matrix pencil
\begin{equation}\label{eq:matpencil}
\sigma \Rtt(\lambda) -  \Ttt(\lambda)\equiv 
\sigma\Stt_{J+1:\infty,1:\infty} \Dtt(\lambda) \left(\Stt_{J+1:\infty,1:\infty}\right)^*  -  \Ttt(\lambda)\,,
\end{equation}
and, if $J>0$, the correspoding eigenvector $\cvect$ of \eqref{eq:matpencil} satisfies the orthogonality condition
\begin{equation}\label{eq:matorth}
\Stt_{1:J,1:\infty}\Dtt(\lambda)\left(\Stt_{J+1:\infty,1:\infty}\right)^*\cvect = \zerovect. 
\end{equation}\label{theorem:10.1}
\end{thm}

\begin{proof} Let $\lambda\in[\kappa_J,\kappa_{J+1})$ be an eigenvalue of Problem (P), let $u(x,\yb)\in H^1_*(\Omega)$
be a corresponding eigenfunction, and let its traces $f$ and $g$ be defined by \eqref{eq:fg}. Then the restriction of $u$ to the 
cylindrical ends $\Ccal$ is a square integrable solution of \eqref{eq:lapl} and therefore $\Pcal_J f = \Pcal_J g = 0$ (or
$\Qcal_J g = g$) and 
$-\Rcal^\Ccal_\lambda g = -\Rcal^\Ccal_\lambda \Qcal_J g = f$, with $\Rcal\Ccal_\lambda$ given by \eqref{eq:rlamC}. Similarly, from the 
restriction on the interior part $\Omega_0$, we obtain that $f =  \Rcal^0_\lambda g$ which immediately leads to 
\begin{equation}\label{eq:ev1}
\Acal_\lambda(-1)g=(-\Rcal_\lambda^\Ccal  - \Qcal_J\Rcal_\lambda^0)g = 0
\end{equation}
and additionally implies the orthogonality
condition 
\begin{equation}\label{eq:orth}
\Pcal_J \Rcal^0_\lambda g = 0\,.
\end{equation}

If we now write $g=\sum_{k=1}^\infty \cvect_k w_{k+J}$ with $\cvect_k = (g,w_{k+J})_\Gamma$ and use the definitions
\eqref{eq:Rklinner}--\eqref{eq:matform}, then by elementary matrix manipulations \eqref{eq:ev1} turns into
$(-\Rtt(\lambda)-\Ttt(\lambda))\cvect = 0$ and \eqref{eq:orth} becomes \eqref{eq:matorth}.
\end{proof}

\section{More on eigenvalues of the operator pencil \eqref{eq:pencildef}}\label{sec:monotone}
Before we proceed to descrtibing in detail the numerical procedure which implements Theorem~\ref{thm:equiv}, we want to establish some more facts
about the eigenvalues $\sigma_j(\lambda)$ of the operator pencil
$\Acal_\lambda(\sigma) = \sigma\Rcal^\Ccal_\lambda-\Qcal_J\Rcal^0_\lambda$ (for $\lambda\in[\kappa_J, \kappa_{J+1})$). We enumerate the 
pencil eigenvalues $\sigma_j(\lambda)$ in decreasing order: $\sigma_1(\lambda)\ge \sigma_2(\lambda)\ge\dots$.

\begin{lem}\label{lem:monotone} The $\sigma$-eigenvalues of $\Acal_\lambda(\sigma)$ are monotone increasing functions of
$\lambda$ on each interval $[\Lambda_1, \Lambda_2]\subset[\kappa_J, \kappa_{J+1})$ which does not contain any eigenvalues $\mu_m$ of the 
homogeneous Neumann problem \eqref{eq:U}.
\end{lem}

\begin{proof} We established in \S\S~ \ref{sec:NDcyl} and \ref{sec:NDint} that each of the operators $\Rcal^\Ccal_\lambda$ and 
$\Qcal_J\Rcal^0_\lambda$ is monotone increasing (in the sense of corrresponding forms) on each interval satisfying the conditions 
of the Lemma. Thus, for each fixed negative value of $\sigma$, the quadratic form of the pencil $\Acal_\lambda(\sigma)$ is monotone decresing in $\lambda$, 
and, as $\Rcal^\Ccal_\lambda>0$, by the inverse function theorem the eigenvalues $\sigma(\lambda)$ increase as functions of $\lambda$.
\end{proof}

We also want to estimate how many pencil eigenvalue curves  $\sigma(\lambda)$ can cross the line $\sigma=-1$ as $\lambda$ increases in a given interval (or, 
equivalently, to estimate \emph{a priori} how many eigenvalues of Problem (P) can be located in a given interval). We need to introduce some additional 
notation. Let $\nu_m$ denote the eigenvalues of the interior homogeneous Dirichlet problem:
\begin{equation}\label{eq:intD}
 -\Delta U = \nu U \;\; \text{in} \;\; \Omega_0, \quad U=0  \;\; \text{on} \;\;
 \Gamma, \quad BU = 0 \;\; \text{on} \;\; \Gamma_0.
\end{equation}
(This problem is similar to \eqref{eq:U} with the only difference being that an auxilliary Neumann condition on $\Gamma$ is replaced by the Dirichlet one.). 

\begin{lem}\label{lem:estimate} Let $\Lambda_2\in(\kappa_J, \kappa_{J+1})\setminus\{\mu_m\}$, and let 
$\Lambda_1>\max\{\kappa_J, \max\{\mu_m:\mu_m<\Lambda_2\}\}$. (In other words, we choose an interval $[\Lambda_1, \Lambda_2]$ in
such a manner that it does not contain any thresholds or any Neumann eigenvalues $\mu_m$.)   
Then the number of eigenvalues of Problem (P) in the interval $[\Lambda_1, \Lambda_2]$
does not exceed 
\begin{equation}\label{eq:Kest}
K:=\#\{m:\mu_m<\Lambda_2\}-\#\{m:\nu_m<\Lambda_2\}\,.
\end{equation}
\end{lem}

\begin{proof} By Lemma \ref{lem:monotone}, the curves $\sigma_j(\lambda)$ are monotone in $\lambda\in[\Lambda_1, \Lambda_2]$.
We wish to estimate from above the total possible number of times these curves intersect the straight line $\sigma = -1$. Obviously, this does not exceed the 
number of curves which are above $-1$ at the right-hand side of the interval  $[\Lambda_1, \Lambda_2]$, 
namely $K'=\#\{j:-1<\sigma_j(\Lambda_2)\}$. We shall now show that $K'\le K$, where $K$ is given by \eqref{eq:Kest}.

Note that any  $\lambda_\star$ for which a particular curve $\sigma_j(\lambda)$ crosses $\sigma=-1$ is at the same time characterised 
by the fact that zero is in the spectrum of the operator 
$\Acal_{\lambda_\star}(-1)\equiv -\Rcal_{\lambda_\star}^\Ccal  - \Qcal_J\Rcal_{\lambda\star}^\Ccal$. 
Thus,
\begin{equation}\label{eq:nlambda}
\begin{split}
K' &\le \#\{\text{positive eigenvalues of }-\Rcal_{\Lambda_2}^\Ccal  - \Qcal_J\Rcal_{\Lambda_2}^0\}\\
&=\#\{\text{negative eigenvalues of }\Rcal_{\Lambda_2}^\Ccal  + \Qcal_J\Rcal_{\Lambda_2}^0\}\,.
\end{split}
\end{equation}

Since $\Rcal_{\lambda}^\Ccal$ is non-negative for any $\lambda$, we can majorize the right-hand side of \eqref{eq:nlambda}
by omitting  it there (and also dropping the projector $\Qcal_J$). Thus,
\begin{equation}\label{eq:nlambda1}
K'\le \#\{\text{negative eigenvalues of }\Rcal_{\Lambda_2}^0\}\,.
\end{equation}
But by \cite {Safarov} (which is in turn a generalisation of \cite{Friedlander}, see also \cite{Agranovich}) there is a known relation between the number of 
negative eigenvalues of an interior Neumann-to-Dirichlet operator, and the differences of the Neumann and Dirichlet counting 
functions:
\begin{equation}
\#\{\text{negative eigenvalues of }\Rcal_{\Lambda_2}^0\}=\#\{m:\mu_m<\Lambda_2\}-\#\{m:\nu_m<\Lambda_2\}=K\,.
\end{equation}
\end{proof}
 
\begin{rem}
The reader may notice that the estimate \eqref{eq:Kest} does not in fact use $\Lambda_1$ and is therefore just 
an estimate on the counting function of the eigenvalues of Problem (P). One can easily improve it if one knowns an estimate from 
below on the number of negative eigenvalues of the operator $\Rcal_{\Lambda_1}^\Ccal  + \Qcal_J\Rcal_{\Lambda_1}^0$, say
\[
\#\{\text{negative eigenvalues of }\Rcal_{\Lambda_1}^\Ccal  + \Qcal_J\Rcal_{\Lambda_1}^0\}\ge K_1
\]
(e.g. obtained numerically). Then in the statement of Lemma \ref{lem:estimate} one can replace $K$ by $K-K_1$. Further improvements
are also possible by introducing the counting functions of a generalised Neumann problem instead of the ordinary Neumann problem, 
but we do not discuss them here. 
\end{rem}

\section{Numerical algorithm}\label{numalg}

Theorem \ref{theorem:10.1} gives the characterisation of the eigenvalues of Problem (P)
which forms the basis of our simple numerical approach to their calculation. We emphasize
that this approach depends only on having software capable of the following:
\begin{itemize}
\item finding eigenvalues and eigenfunctions for problems on bounded domains with mixed 
 boundary conditions (e.g. MATLAB PDE Toolbox, FEMLAB);
\item performing quadratures on boundaries and cross-sections of domains;
\item finding the first few eigenvalues of a matrix pencil;
\item finding the zeros of a monotone real-valued function of a real variable.
\end{itemize}
Our algorithm in its basic form is as follows:
\begin{algorithm}\label{algorithm:1}
\

We remind the reader that the $\mu_m$ are the eigenvalues of our PDE on $\Omega_0$ with Neumann
boundary conditions on the interface $\Gamma$ (problem (\ref{eq:intN}); see also Fig. \ref{fig:domain}); 
the $U_m$ are the corresponding eigenfunctions; and the $w_k$ are the eigenfunctions of the 
interface problem \eqref{eq:trans} with corresponding eigenvalues $\kappa_k$. In two dimensions, 
the interface problem \eqref{eq:trans} becomes an ODE problem and the $w_k$ will be explicit 
trigonometric or hyperbolic functions.

\begin{description}
\item[1.]
Choose a sufficiently large $N\in \Nbb$ and
calculate numerical approximations to the Neumann eigenvalues $(\mu_m)_{m=1}^N$ and their 
eigenfunctions $(U_m)_{m=1}^N$ for the inner Neumann problem \eqref{eq:intN}, using an appropriate
method -- e.g. a standard finite element package.

\item[2.] Choose $M\in\Nbb$ such that the traces of the eigenfunctions $(U_m)_{m=1}^N$ 
 on the interface $\Gamma$ are well approximated by linear combinations of the first
 $M$ eigenfunctions $(w_j)_{j=1}^M$ of the interface problem \eqref{eq:trans}.  
 Since the convergence of such expansions is usually exponential it is rare to need 
 $M>20$.

\item[3.] Calculate the first $M$ eigenvalues $(\kappa_j)_{j=1}^M$ and eigenfunctions
 $(w_j)_{j=1}^M$ of the interface problem \eqref{eq:trans} --- either analytically, in terms
 of sines and cosines, as in the case of the two-dimensional examples treated here, or 
 numerically in the case of general higher-dimensional examples. 

\item[4.] Calculate the entries
\[ S_{km} = \left(w_k,U_m|_\Gamma\right)_\Gamma =\int_\Gamma w_k\left.\overline{U_m}\right|_\Gamma;
\]
for $k=1,\ldots,M$, $m=1,\ldots,N$ of the matrix $\Stt$ (see \eqref{eq:Sdef}):
these are the inner products of the interface eigenfunctions $w_k$ with the traces on the interface
of the eigenfunctions $U_m$ of problem (\ref{eq:intN}).

\item[5.] Fix $J$, $0\leq J \ll M$ -- we shall look for eigenvalues $\lambda\in [\kappa_{J},\kappa_{J+1})$.
 Recall that we are using the convention $\kappa_0 = -\infty$, so for eigenvalues below the essential
 spectrum we choose $J=0$. No calculations are carried out at this step.

\item[6.] Calculate the first $N$ Dirichlet eigenvalues of the problem \eqref{eq:intD} 
 and use the estimates in \S~\ref{sec:monotone} to obtain an upper bound $K$ on
 the number of eigenvalues of Problem (P) in $[\kappa_J,\kappa_{J+1})$. Note:
 $K\ll M$.

\item[7.] Define the functions giving the elements $D_{kk}$, $k=1,\dots,N$, 
of the diagonal matrix $\Dtt(\lambda)$ and the elements $T_{kk}$, $k=1,\dots,M$, 
of the diagonal matrix $\Ttt(\lambda)$, by
\[ D_{kk}(\lambda) = (\mu_k-\lambda)^{-1}; \;\;\;
   T_{kk}(\lambda) = (\kappa_k-\lambda)^{-1/2}. \]
Here the $\kappa_k$ are the eigenvalues belonging to the interface eigenfunctions $w_k$, calculated
in Step 3.  No calculations are carried out at this step.

\item[8.] Define a procedure which takes a value $\lambda\in [\kappa_J,\kappa_{J+1})$
and returns the largest $K$ $\sigma$-eigenvalues
\[ \sigma_1(\lambda)\geq\sigma_2(\lambda)\geq \cdots \geq \sigma_{K}(\lambda) \] 
(out of the total of $M-J$ $\sigma$-eigenvalues) of the matrix pencil
\begin{equation}\label{eq:approxpencil}
 \sigma \Stt_{J+1:M,1:N}\Dtt_{1:N,1:N}(\lambda)\left(\Stt_{J+1:M,1:N}\right)^*-\Ttt_{J+1:M,J+1:M}(\lambda).
\end{equation}
$K$ is given in Step 6. No calculations are carried out at this step.

\item[9.] Using the procedure of Step 8 and an appropriate numerical rootfinding algorithm,
solve (if possible) the $K$ nonlinear equations
\begin{equation}\label{eq:sigmaj}
\sigma_j(\lambda) = -1, \;\;\; j = 1,\ldots, K.
\end{equation} 
 Note that care has to be taken if any of the interior Neumann eigenvalues $\mu_m$ 
 lies in $[\kappa_J,\kappa_{J+1})$. Away from the $\mu_m$, the $\sigma_j$ are monotone
 functions of $\lambda$ and so the rootfinding is, in principle, easy.

 Denote by $\lambda_{j,J}$ the solution of \eqref{eq:sigmaj} and by $\cvect_{j,J}$
 the eigenvector of the pencil \eqref{eq:approxpencil} with $\lambda=\lambda_{j,J}$, 
 $\sigma=-1$.

\item[10.] If $J>0$, check the orthogonality conditions 
\begin{equation}
 \Stt_{1:J,1:N}\Dtt_{1:N,1:N}(\lambda)\left(S_{J+1:M,1:N}\right)^*\cvect_{j,J}=\zerovect. 
\label{eq:ortho2}
\end{equation}
(Of course, these can never be checked exactly: at most one can check that
they are `almost' satisfied and that there is therefore either an embedded eigenvalue
or a point of spectral concentration near $\lambda_{j,J}$.)

\end{description} \hfill $\Box$
\end{algorithm}

In the absence of error from the numerical approximation of the eigenvalues and
eigenfunctions, the accuracy of this algorithm is determined by the choice of $N$ 
at the very first step. (The choice of $M$ in step 2 is generally less problematic
as the transversal eigenfunctions provide rapidly converging approximations to
smooth functions on the interface.) $N$ is the number of terms at which the sum
in the expression (\ref{eq:mainexpr}) for $R_{k\ell}(\lambda)$ is truncated, and
unfortunately this series can be rather slowly convergent. We accelerate the 
convergence with a cheap trick. Suppose that for some fixed $\lambda_0$ we already 
know $R_{k\ell}(\lambda_0)$. Then (\ref{eq:mainexpr}) yields
\begin{equation}
R_{k\ell}(\lambda)-R_{k\ell}(\lambda_0) = 
 \sum_{m=1}^\infty
\frac{\lambda-\lambda_0}{(\mu_m-\lambda)(\mu_m-\lambda_0)}(\phi_k,U_m|_\Gamma)_\Gamma\cdot(U_m|_\Gamma,\phi_\ell)_\Gamma\,.
\label{eq:imprexpr}
\end{equation}
The series in (\ref{eq:imprexpr}) converges much more rapidly than the series in
(\ref{eq:mainexpr}). Moreover the numbers $R_{k\ell}(\lambda_0)$ are easily
calculated from (\ref{eq:startexp}) (with $\lambda=\lambda_0$) by solving the 
boundary value problems (\ref{eq:intPhi}) numerically (for $\lambda=\lambda_0$) 
using any reasonable discretization scheme for elliptic PDEs. Our Main 
Algorithm is then modified by changing steps 8 and 10 as follows:
\begin{description}

\item[8'.] Let $\Rtt^0$ denote the matrix with entries $R_{k\ell}(\lambda_0)$ defined
in (\ref{eq:Rklinner}) (with $\lambda$ replaced by $\lambda_0$).
Then the matrix pencil (\ref{eq:approxpencil}) is replaced by
\begin{equation}
\begin{array}{l}
\sigma \left[ \Rtt^0_{J+1:M,J+1:M}+\Stt_{J+1:M,1:N}\left\{\Dtt_{1:N,1:N}(\lambda)-\Dtt_{1:N,1:N}(\lambda_0)\right\}\left(\Stt_{J+1:M,1:N}\right)^*\right] \\
 \hspace{5cm} -\Ttt_{J+1:M,J+1:M}(\lambda). \end{array}
\label{eq:approxpencil2}
\end{equation}

\item[10'.] The orthogonality condition (\ref{eq:ortho2}) is replaced by
\begin{equation}
 \left[\Rtt^0_{1-J:0,1:N}+\Stt_{1:J,1:N}\left\{\Dtt_{1:N,1:N}(\lambda)
-\Dtt_{1:N,1:N}(\lambda_0)\right\}\left(S_{J+1:M,1:N}\right)^*
\right]\cvect_{j,J}=\zerovect. 
\label{eq:ortho3}
\end{equation}
Note that $\Rtt^0_{1-J:0,1:N}$ makes sense: in view of our choice of basis functions
$\phi_k = w_{k+J}$, its $(k,j)$ element is the numerical approximation
to $(\Rcal_{\lambda_0}^0w_k,w_{j+J})_\Gamma$.
\end{description}

The convergence acceleration trick can be repeated if one is prepared to evaluate
the matrix elements $R_{k\ell}$ numerically for several different $\lambda$. However
for the numerical examples which we examined here, one application was always
sufficient.

\subsection{Example 1: Bent waveguide}

\begin{figure}[htb!]
\begin{center}
%\resizebox{425pt}{250pt}{\includegraphics*{twistedtube}}
\resizebox{250pt}{!}{\includegraphics*{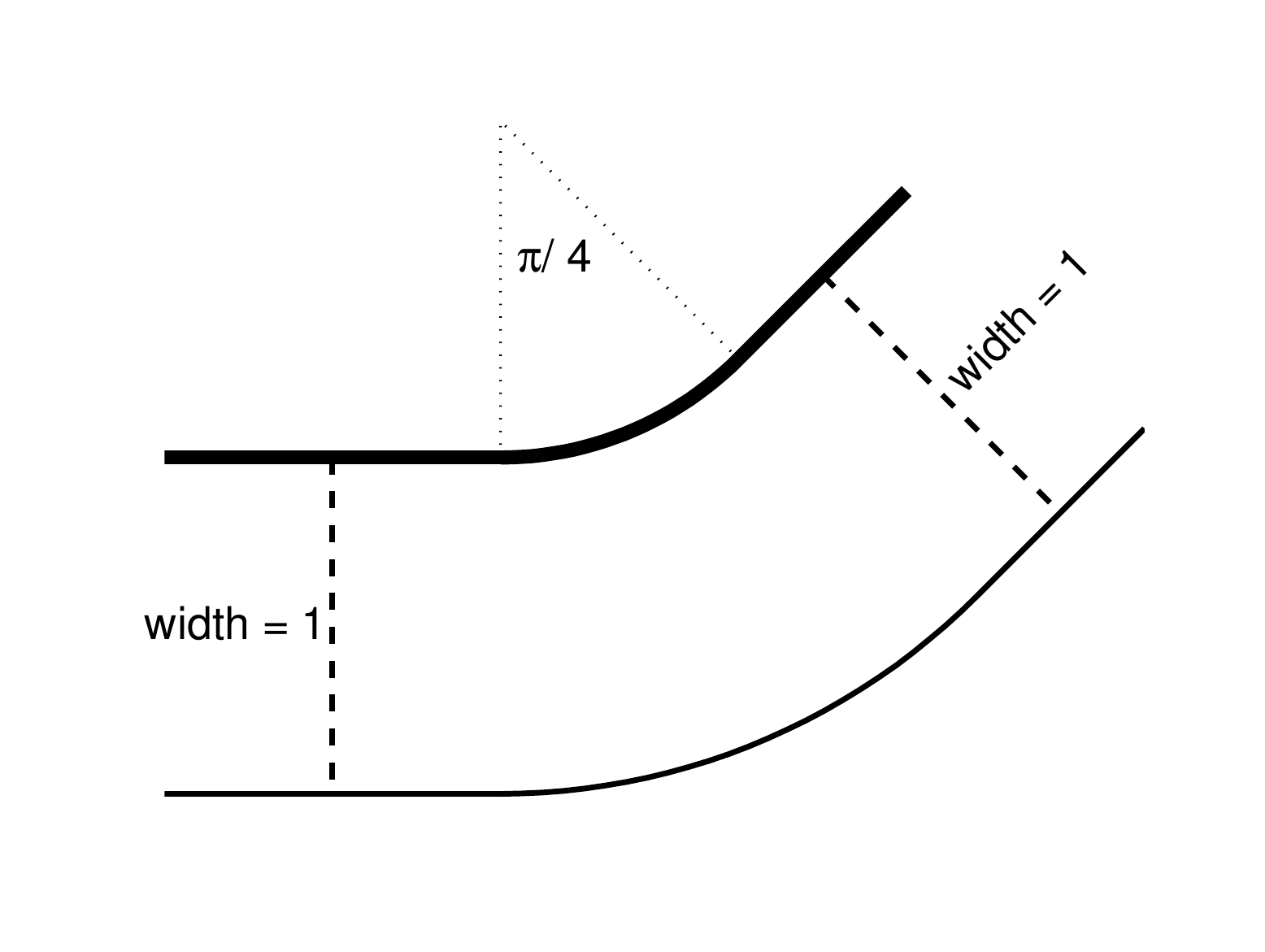}}
\caption{Bent Waveguide\label{fig:mm1}}
\end{center}
\end{figure}

We consider the Laplacian in a waveguide in $\Rbb[2]$ of width 1, 
bent through an angle $\pi/4$ as indicated in Fig. \ref{fig:mm1}
(the inner radius of curvature of the bend is 1, the outer
radius 2). On the side of the waveguide with the smaller
radius of curvature on the bend, Dirichlet conditions are
imposed; on the other side the boundary conditions are of
Neumann type. 

We know that the essential spectrum is determined
by the spectrum of the problem on the cylindrical ends only, and is
therefore the same as for the unbent waveguide:
\[ \sigma_{\mathrm{ess}} = [\pi^2/4,+\infty). \]
It may also be shown that for the unbent waveguide the spectrum is
purely absolutely continuous. However this changes when the waveguide
is bent. Krej\v{c}i\v{r}\'{\i}k and K\v{r}\'{\i}\v{z} \cite{KK} show that bending the
waveguide in the direction of the Dirichlet boundary condition will
cause an eigenvalue to appear below the essential spectrum. Using
the method described in \S~\ref{numalg}, this eigenvalue can
be located quite precisely. For this purpose the inner domain 
$\Omega_0$ is taken to be the part of the domain between the two
dashed lines in Fig. \ref{fig:mm1}. The basis on the cut-off
boundary $\Gamma$ of $\Omega_0$ (the part of $\partial\Omega_0$
composed of the two dashed lines in Fig. \ref{fig:mm1}) is defined as 
follows. Let $s$ be the local transverse coordinate on $\Gamma$ measured 
with $s=0$ on the `Dirichlet' side.
We denote by $\phi_{2j+1}$, $j=0,1,2,\ldots$, the basis functions
supported on the vertical portion of $\Gamma$, and
by  $\phi_{2j}$, $j=1,2,\ldots$, the basis supported on the part of 
$\Gamma$ which is at angle $\pi/4$ to the horizontal. 
In particular, $\phi_{2j}$ is zero on the vertical part
of $\Gamma$ and $\phi_{2j+1}$ is zero on the part of $\Gamma$ which
is at angle $\pi/4$ to the horizontal. On their supports, the
$\phi_j$ are given by
\[ \phi_{j}(s) = \left\{ \begin{array}{lr} 
                            \sin(j\pi s/2) & \mbox{($j$ odd)}, \\
                            \sin((j-1)\pi s/2) &  \mbox{($j$ even)}. 
                         \end{array}\right. \]
Table \ref{table:mm2} shows the results of  computations at four
different levels of accuracy.

\begin{table}[tbh]
\begin{center}\normalsize\small
\begin{tabular}{lllc}\hline
\multicolumn{3}{c}{accuracy} & \\
\multicolumn{3}{c}{$\overbrace{\qquad\qquad\qquad\qquad\qquad\qquad\qquad\qquad\qquad\qquad\qquad}$}&eigenvalue found \\ \hline
silly&one mesh refinement&sum over $\lambda_j\leq 10$ in (\ref{eq:mainexpr}) & none found \\ 
low&three mesh refinements&sum over $\lambda_j\leq 50$ in (\ref{eq:mainexpr})& 2.3461  \\
medium&four mesh refinements&sum over $\lambda_j\leq 100$ in (\ref{eq:mainexpr})& 2.3459  \\
high&five mesh refinements&sum over $\lambda_j\leq 200$ in (\ref{eq:mainexpr})& 2.3454 \\ \hline
\end{tabular}
\end{center}
\caption{Levels of accuracy and eigenvalue found below the essential spectrum, caused by bending the waveguide.
\label{table:mm2}}
\end{table}

%\begin{table}[tbh]
%\hspace{3cm}
%\begin{tabular}{|l|l|}\hline
%Silly & 1 mesh refinement; sum over $\lambda_j\leq 10$ in (\ref{eq:mainexpr}); \\ \hline
%Low & 3 mesh refinements; sum over $\lambda_j\leq 50$ in (\ref{eq:mainexpr}); 
%\\ \hline
%Medium & 4 mesh refinements; sum over $\lambda_j \leq 100$ in (\ref{eq:mainexpr}); 
%\\ \hline
%High & 5 mesh refinements; sum over $\lambda_j \leq 200$ in (\ref{eq:mainexpr}).
%\\ \hline
%\end{tabular} 
%\caption{Levels of Accuracy\label{table:acc}}
%\end{table}

\section{Generalization to resonances}
\label{sec:res1}
The algorithm in \S~\ref{numalg} can be generalized to the calculation of 
resonances. There are two steps in this generalization.

Firstly, resonances are generally complex. (When a resonance becomes real,
it usually becomes either an imbedded eigenvalue, which we have already discussed,
or a point of spectral concentration --- a local maximum of the derivative of one of the 
spectral measures associated with the operator.)
Since resonances are complex, the condition $\lambda\in [\kappa_J,\kappa_{J+1})$
introduced in \S~\ref{sec:NDcyl} for the calculation of imbedded eigenvalues
no longer makes sense. Of course the equation $-\Delta V = \lambda V$ can still
be solved in the cylindrical ends, following the procedure in \S~\ref{sec:NDcyl},
and the solution expressed in the form 
\[ 
V(x,\yb) = \sum_{j=1}^\infty c_j(x) w_j(\yb). 
\]
The functions $c_j$ are now chosen to be
\[
c_j(x)=c_j^+(x)=-\frac{(g,w_j)}{\sqrt{\kappa_j-\lambda}}\exp\left((x-1)\sqrt{\kappa_j-\lambda}\right)\,.
\]
(cf. \eqref{eq:cjpm}). The reason for such a choice is that for calculation of resonances we drop the condition 
$V\in H^1_*(\Ccal)$ and we explicitly
choose the branch of the square root such that $c_j\not\in L_2(0,\infty)$. Thus, if we seek
resonances (for definiteness, in the lower half-plane $\Im\lambda<0$), we choose the branch of the square root such that
$\Re(\sqrt{\kappa_j-\lambda})>0$.

The expression \eqref{eq:rlamC} remains formally unchanged: however one must choose the
branch of the square root as above. A similar remark applies to the
matrix $\Ttt$ defined in \eqref{eq:Tdef}. With this proviso the matrix pencil 
\eqref{eq:approxpencil2} is formally unchanged, though $J$ must be assigned the value
$0$. 

With these modifications, the resonance approximations are the values of $\lambda$ such 
that some $\sigma$-eigenvalue $\sigma_j(\lambda)$ of the pencil (\ref{eq:approxpencil2}) 
satisfies $\sigma_j(\lambda)=1$. However calculating the resonances in this way is
impractical. Even if the $\sigma_j$ are analytic functions of $\lambda$, which is not
guaranteed, it is difficult to order the numerically calculated values of the $\sigma_j$
to correspond to an analytic ordering.

The simple approach which we used here was to calculate the condition number
of the pencil 
\begin{equation}
\left[ \Rtt^0_{1:M,1:M}+\Stt_{1:M,1:N}\left\{\Dtt_{1:N,1:N}(\lambda)-\Dtt_{1:N,1:N}(\lambda_0)\right\}\left(\Stt_{1:M,1:N}\right)^*\right]-\Ttt_{1:M,1:M}(\lambda).
\label{eq:approxpencil3}
\end{equation}
as a function of $\lambda$ and make a contour plot. Local maxima are `suspicious
points' as the condition number will be infinite at a resonance. (The converse,
however, is not true: not all local maxima are resonances.) Once suspicious points
are located approximately one can `zoom in' and examine them in more detail.

\subsection{Example 2: obstructed waveguide}
\begin{figure}
\begin{center}
\resizebox{300pt}{130pt}{\includegraphics*{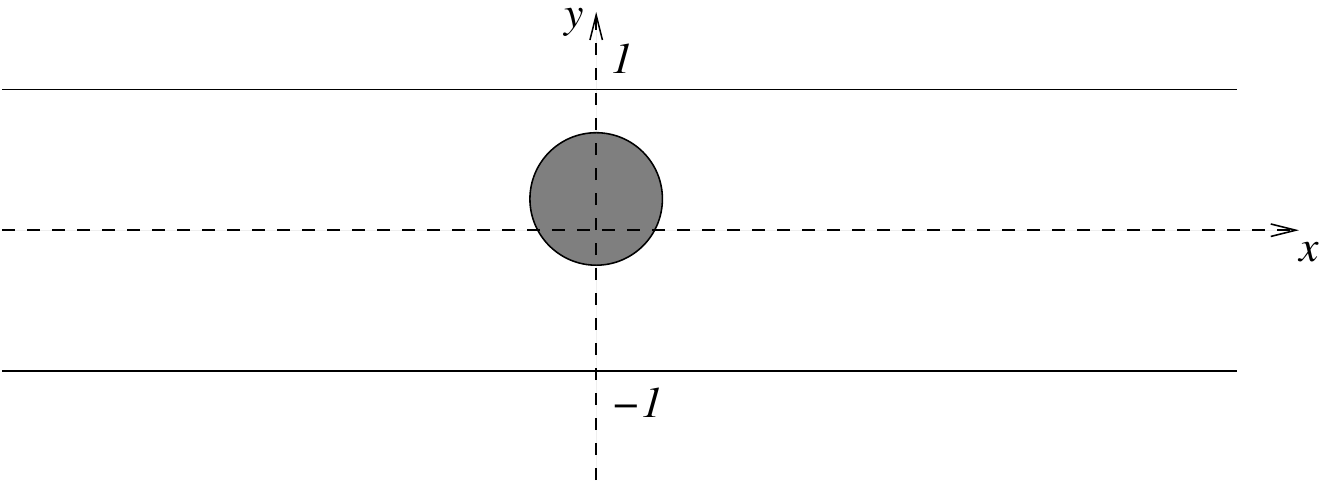}}
\end{center}
\caption{Waveguide obstructed by an obstacle\label{fig:waveguide}}
\end{figure}

This problem is considered in detail by Aslanyan, Parnovski and
Vassiliev \cite{APV}: a waveguide in $\Rbb[2]$ in the strip
$|y|<1$, obstructed by a plane obstacle ${\cal O}$ --- see Fig. \ref{fig:waveguide}. 
For the particular experiments here, we consider an obstacle which is symmetric about the 
$y$-axis. The domain can then be reduced to the strip
\[ \{ (x,y)\in\Rbb[2]\; : \; x>0, \; |y|<1 \} \setminus \overline{\cal O}. \]
Following Aslanyan et al. we take as the obstacle a disc of centre
$(0,\delta)$ and radius $R$, with $|\delta|+R<1$:
\[ {\cal O} = \{ (x,y) \; : \; x^2 + (y-\delta)^2 < R^2\}. \]
When $\delta=0$ the problem has an eigenvalue imbedded in the essential
spectrum. When $\delta$ is perturbed from zero, the eigenvalue evolves 
into a resonance, analytically as a function of $\delta$ (see \cite{APV}). 
For the numerical experiments we truncated the strip at $x=1$:
thus, we took as our domain $\Omega_0$ the set
\[ \{ (x,y)\in\Rbb[2]\; : \; 0<x<1, \; |y|<1 \} \setminus \overline{\cal O}. \]
Table \ref{table:mm1} shows the results for various values of the
numerical experiment parameters. Each is quoted for four
`levels of accuracy' defined in Table \ref{table:mm2}.

\begin{table}[tbh!]\normalsize\small
\begin{center}\begin{tabular}{clclc}
\hline
 & \multicolumn{2}{c}{$R=0.3$}& \multicolumn{2}{c}{$R=0.5$}\\
 & \multicolumn{2}{c}{$\overbrace{\qquad\qquad\qquad\qquad\qquad\qquad}$}& \multicolumn{2}{c}{$\overbrace{\qquad\qquad\qquad\qquad\qquad\qquad}$}\\
$\delta$ &  eigenvalue & accuracy & eigenvalue & accuracy \\
\hline
    & 1.50866 & S & 1.39806  & S  \\
    & 1.50499 & L & 1.39157 & L  \\
0.0 & 1.50489 & M & 1.39139 & M \\
    & 1.50486 & H   & 1.39134 & M \\
    & 1.5048 & A & 1.3913 & A  \\
    \multicolumn{5}{c}{\ }\\
    & $1.5120+10^{-4}\mathrm{i}$ & S   &  $1.4043+8\times 10^{-4}\mathrm{i}$ & S    \\
    & $1.5080+10^{-4}\mathrm{i}$ & L  & $1.3981+9\times 10^{-4}\mathrm{i}$ & L  \\
0.1 & $1.5079+10^{-4}\mathrm{i}$ & M & $1.3979+9\times 10^{-4}\mathrm{i}$ & M  \\
    & $1.5078+10^{-4}\mathrm{i}$ & H  & $1.3979+9\times 10^{-4}\mathrm{i}$ & H   \\ 
    & $1.5102+10^{-4}\mathrm{i}$ & A & $1.3998+9\times 10^{-4}\mathrm{i}$ & A  \\ 
    \multicolumn{5}{c}{\ }\\
    & $1.5203+4\times 10^{-4}\mathrm{i}$ & S & $1.4236+3.46\times 10^{-3}\mathrm{i}$ & S \\
    & $1.5167+5\times 10^{-4}\mathrm{i}$ & L & $1.4180+3.89\times 10^{-3}\mathrm{i}$ & L \\ 
0.2 & $1.5165+5\times 10^{-4}\mathrm{i}$ & M & $1.4178+3.89\times 10^{-3}\mathrm{i}$ & M \\ 
    & $1.5165+5\times 10^{-4}\mathrm{i}$ & H  & $1.4178+3.90\times 10^{-3}\mathrm{i}$ & H  \\ 
    & $1.5188+5\times 10^{-4}\mathrm{i}$ & A & $1.4196+3.93\times 10^{-3}\mathrm{i}$ & A \\ \hline
\end{tabular}
\end{center}
\caption{Experiments on the obstructed waveguide. (Accuracy: S, silly; L, low; M, medium; H, high; A,  from Aslanyan et al. \cite{APV}.)\label{table:mm1}}
\end{table}

The agreement with the results of Aslanyan et al. is good for the case of 
obstacle displacement $\delta = 0$ (the case of an embedded eigenvalue). In the
cases $\delta>0$ the agreement is less good, but the error is no worse than $2\times 10^{-3}$.
The values which we calculate change very little between the medium and high accuracy
cases.

\section{Generalizations --- resonances of a scatterer}
\label{sec:res2}
\begin{figure}[htb!]
\begin{center}
\resizebox{200pt}{200pt}{\includegraphics*{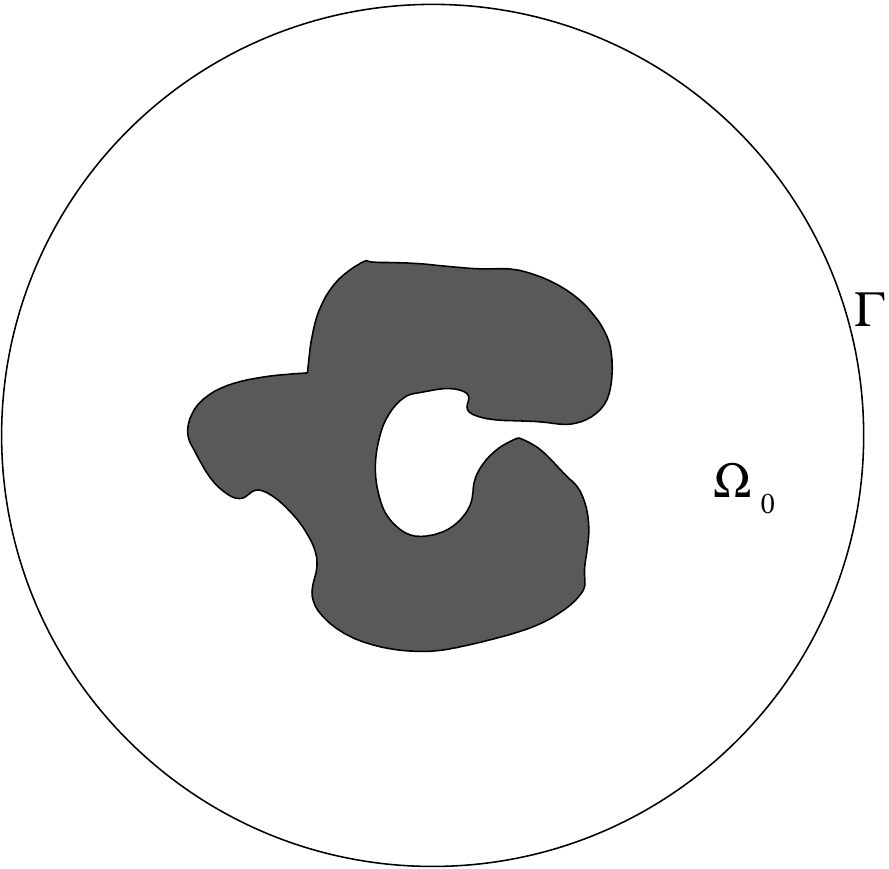}}
\caption{Scattering by an obstacle. The transverse boundary $\Gamma$ is replaced
by a circle of radius $R$, the cylindrical ends by the region $x^2+y^2>R^2$.\label{fig:mm3}}
\end{center}
\end{figure}

In\S~\ref{sec:res1} we described how to calculate resonances in the
case when the domain has cylindrical ends. In fact a similar approach can
be used to calculate resonances due to scattering by an obstacle. We consider
a situation similar to that in Fig. \ref{fig:mm3}, in which the domain is
the region exterior to an obstacle (the shaded region). The cylindrical ends
are now replaced by the region $x^2+y^2>R^2$ and the boundary $\Gamma$ is the
circle of centre $(0,0)$, radius $R$. 

As in \S~\ref{sec:res1}, we can compute resonances by replacing the
map $\Rcal^\Ccal_\lambda$ of \eqref{eq:rlamC} by the appropriate Neumann
to Dirichlet map for the exterior region. A simple separation of variables
shows that, in terms of appropriate Hankel functions $H_n$,
\begin{equation}\label{eq:rlamO}
\Rcal^\Ccal_\lambda g = \frac{1}{\sqrt{2\pi}}\sum_{n\in{\mathbb Z}}g_n \frac{H_n(R\sqrt{\lambda})}{
 \sqrt{\lambda}H_n'(R\sqrt{\lambda})}\mbox{e}^{in\theta}, 
\end{equation}
in which the Fourier coefficients $g_n$ are
\[
g_n = \frac{1}{\sqrt{2\pi}}\int_0^{2\pi}g(\phi)\mbox{e}^{-in\phi}{\rm d}\phi. 
\]
Suppose that we are interested in resonances which lie in the half-plane
$\Im(\lambda)<0$. Then choosing the branch of the square root so that
$\Im(\sqrt{\lambda})<0$, the required Hankel functions are the ones which 
grow as
\begin{equation}
\label{eq:hankel}
 H_n(r\sqrt{\lambda}) \sim 
 r^{-1/2}\mbox{e}^{-r\im\sqrt{\lambda}},
\end{equation}
as $r\rightarrow +\infty$. These
are the Hankel functions of the first kind. 

The matrix $\Ttt$ of \eqref{eq:Tdef} is now replaced by a matrix which 
is clearly still diagonal. If we order the functions $\{ \mbox{e}^{in\theta} \, | \, n\in 
{\mathbb Z}\} $ in the order 
\[ \{ 1, \mbox{e}^{i\theta}, \mbox{e}^{-i\theta},
 \mbox{e}^{2i\theta}, \mbox{e}^{-2i\theta}, \ldots \} \]
then the corresponding matrix $\Ttt$ is
\begin{equation}
\label{eq:Tdef2}
\Ttt = \mbox{diag}\left( \frac{H_0(R\sqrt{\lambda})}{\sqrt{\lambda}H_0'(R\sqrt{\lambda})},
\frac{H_1(R\sqrt{\lambda})}{\sqrt{\lambda}H_1'(R\sqrt{\lambda})},
\frac{H_{-1}(R\sqrt{\lambda})}{\sqrt{\lambda}H_{-1}'(R\sqrt{\lambda})},
\frac{H_2(R\sqrt{\lambda})}{\sqrt{\lambda}H_2'(R\sqrt{\lambda})},
\frac{H_{-2}(R\sqrt{\lambda})}{\sqrt{\lambda}H_{-2}'(R\sqrt{\lambda})},\ldots \right) \, .
\end{equation}
The resonances can then be calculated following the same procedure as in \S~\ref{sec:res1}.

\subsection{Example 3: a cavity resonance problem}
\begin{figure}[htb!]
\begin{center}
\resizebox{270pt}{230pt}{\includegraphics*{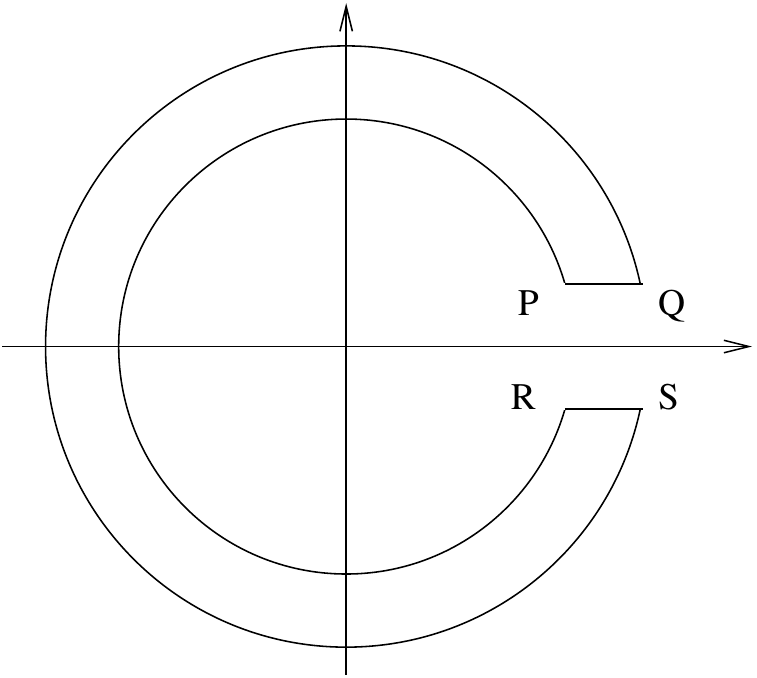}}
\caption{Cavity Resonator. The points $P$, $Q$, $R$ and $S$ have
coordinates $(1,\epsilon)$, $(1.1,\epsilon)$, $(1.1,-\epsilon)$ and $(1,-\epsilon)$
respectively, and are joined by two circular arcs\label{fig:mm2}}
\end{center}
\end{figure}

We consider scattering by a C-shaped barrier almost enclosing a cavity
as shown in Fig. \ref{fig:mm2}. 

On the boundary of the domain we impose Dirichlet boundary
conditions. For the Laplacian on this domain,
$\sigma_{ess}=[0,+\infty)$. By simple separation of variables
it is easy to show that no non-trivial solution of the Poisson 
equation $-\Delta u = \lambda u$ on this exterior domain can lie 
in $L_2$ for $\lambda\geq 0$, no matter what boundary conditions 
it may satisfy. Thus there can be no embedded eigenvalues.

Nevertheless it is known (see \cite{BrHiMa}) that
for small $\epsilon$ this problem possesses a resonance close 
to each Dirichlet eigenvalue for the Laplacian on the unit disc. 
We wish to compute these resonances.

We first create an `inner domain' by introducing an artificial
outer boundary which is a circle centred at the origin with
radius $R>1.1$. It is important that $R$ not be too large:
as $R\rightarrow\infty$, the (Neumann and Dirichlet) eigenvalues 
of the inner domain become dense in any compact subset of 
$[0,\infty)$, and so in particular it becomes very difficult
to distinguish resonances close to the real axis from Neumann
eigenvalues of the inner domain. For numerical reasons it is's
also important that $R$ not be too close to $1.1$: otherwise
the annular region $1.1 < \sqrt{x^2+y^2} < R$ becomes very narrow,
requiring a very fine mesh, and the numerical calculations 
require too much computing time. We carried out our calculations
for $R=1.3$ and, as a check, for $R=1.5$.

%Next, we require the ND map for the exterior domain $\sqrt{x^2+y^2}>R$.
%For $\Im(\lambda) < 0$ this can be written explicitly in terms
%of Hankel functions: 
%\[ \Lambda_{ext}(\lambda)\sum_{n\in {\Bbb Z}}f_n\mbox{e}^{in\theta}
% = \sum_{n\in {\Bbb Z}}f_n\frac{H_n(R\sqrt{\lambda})}{\sqrt{\lambda}H_n'(R\sqrt{\lambda})}
% \mbox{e}^{in\theta}, \]
%in which $H_n$ is the unique (up to scalar multiples) Hankel function
%such that $H_n(z)\rightarrow 0$ as $\Im(z)\rightarrow -\infty$, and
%$\sqrt{\lambda}$ denotes the unique square root of $\lambda$ such that
%$\Im(\sqrt{\lambda})<0$ when $\Im(\lambda)<0$. It is therefore easy
%to calculate the matrix elements of $\Lambda_{ext}$ with respect to
%the orthonormal basis $\displaystyle \left\{ (2\pi R)^{-1/2}\mbox{e}^{in\theta} \, | \, 
%n \in {\Bbb Z} \right\} $, and of course the matrix representation
%of $\Lambda_{ext}$ turns out to be diagonal. 
%
%In order to find `first channel' resonances one replaces the Hankel
%functions for $n=1$ by the corresponding exponentially growing Hankel
%function (uniquely d%etemined by exponential decay exponentially when $\Im(z)\rightarrow
%+\infty$). Resonances in higher channels are calculated by the same
%operation for higher values of $n$. We calculated only first channel
%resonances here.

\begin{figure}[htb!]
\begin{center}
\resizebox{340pt}{300pt}{\includegraphics*{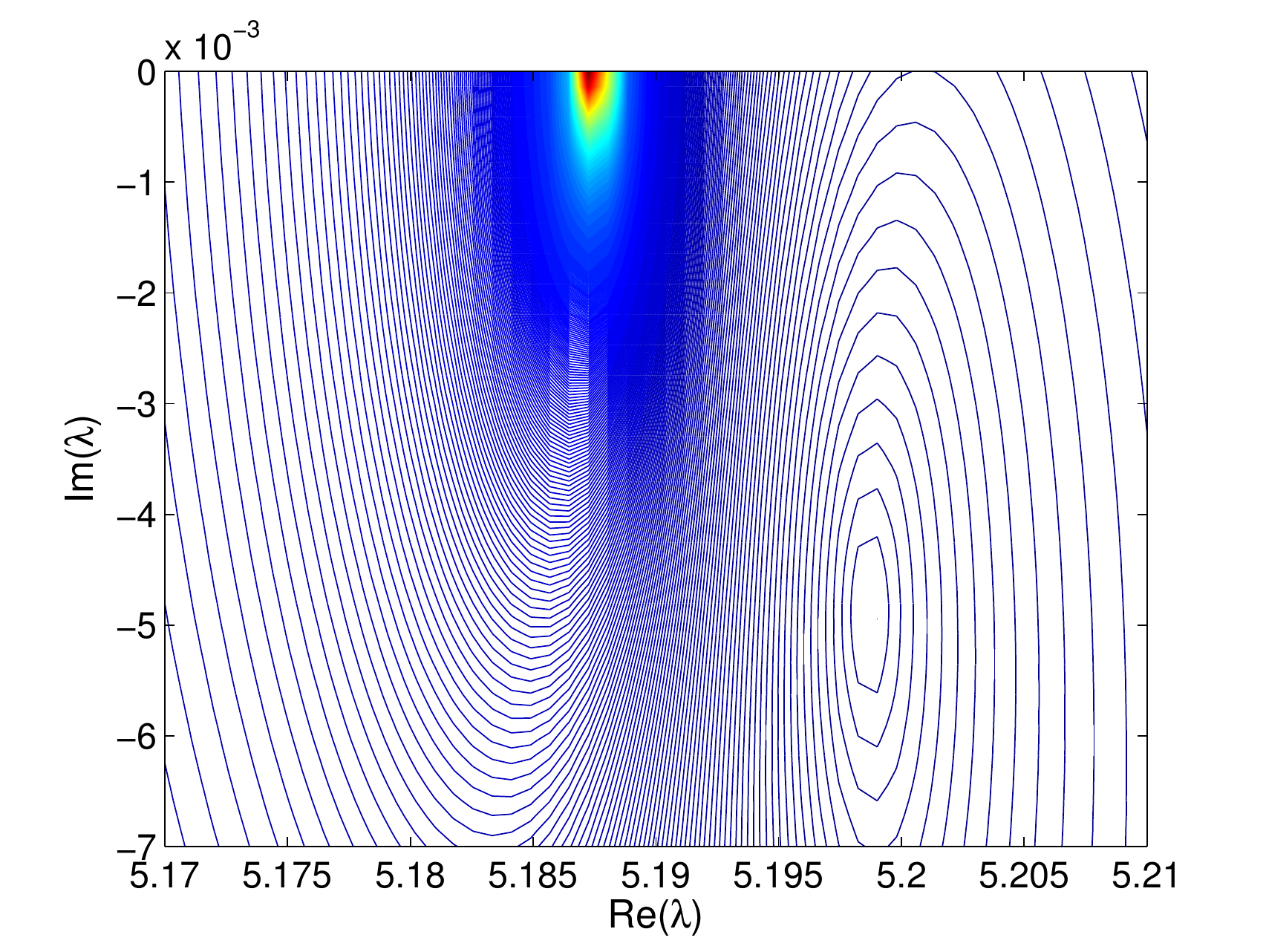}}
\caption{Resonance very close to a Neumann eigenvalue for the interior domain;
$R=1.5$, $\epsilon=0.3$. At lower resolutions these are virtually indistinguishable. 
This is a problem, as not all Neumann eigenvalues are close to resonances! The contours
are the level sets of a determinant which is zero at the resonance $5.199-5\times 10^{-3}\mathrm{i}$. 
Unfortunately the same determinant also has a pole at the nearby Neumann eigenvalue 
$5.187$.\label{fig:mm3.14}}
\end{center}
\end{figure}

The Dirichlet problem on the unit disc has eigenvalues 5.783 (simple)
and 14.682 (multiplicity 2). We started by calculating the first channel
resonances of the C-shape resonator closest to these eigenvalues. The results
are shown in Table \ref{table:mm3}. 
\begin{table}[tb!]
\begin{center}\small
\begin{tabular}{ccc}\hline
                 & $\epsilon=0.2$ & $\epsilon = 0.3$ \\ \hline
 medium & $5.5158-4.8\times 10^{-4}\mathrm{i}$ & $5.1987-0.0049\mathrm{i}$ \\
                 &  $13.8574-4\times 10^{-4}\mathrm{i}$ & $12.7505-0.1045\mathrm{i}$ \\ \hline
high & $5.5141-4.6\times 10^{-4}\mathrm{i}$ & $5.1976-0.0051\mathrm{i}$ \\
               &  $13.8554-9.7\times 10^{-3}\mathrm{i}$ &  $12.8165-0.1025\mathrm{i}$  \\ \hline
\end{tabular}
\end{center}
\caption{Resonance calculations for $C$-shaped resonator. (Levels
of accuracy: `medium' means 4 mesh refinements, $R=1.5$ and $\lambda_{\text{max}}=50$;
`high' means 5 mesh refinements, $R=1.3$ and $\lambda_{\text{max}}=100$.)\label{table:mm3}}
\end{table}
The resonance closest to the first of these eigenvalues is rather difficult to calculate accurately
because it is so close to a pole of $\Rcal^0_\lambda$, the Neumann-to-Dirichlet
map of the interior domain whose poles are the Neumann eigenvalues
of the interior domain. See Fig. \ref{fig:mm3.14}, which shows a contour plot
of 
\[ \left|\det(L_{\text{inner}}(\lambda)L_{\text{outer}}^{-1}(\lambda)-I)\right|, \]
where, in view of (\ref{eq:approxpencil3}), $L_{\text{inner}}$ is given by
\[ L_{\text{inner}}(\lambda) = \Rtt^0_{1:M,1:M}+\Stt_{1:M,1:N}
 \left\{\Dtt_{1:N,1:N}(\lambda)-\Dtt_{1:N,1:N}(\lambda_0)\right\}\left(\Stt_{1:M,1:N}\right)^*,
\]
while $L_{\text{outer}}(\lambda)$ is the analytic continuation of $\Ttt_{1:M,1:M}(\lambda)$
already described. The point close to $5.2-5\times 10^{-3}\mathrm{i}$ is the zero of
the determinant, and is the approximate resonance. It is not easy to see this
zero so close to the nearby pole. This requires sampling of the determinant on
a very fine grid, which is expensive.

The resonance close to the second eigenvalue is not so close as to cause the same
numerical difficulties. 

\begin{figure}[htb!]
\begin{center}
\resizebox{250pt}{220pt}{\includegraphics*{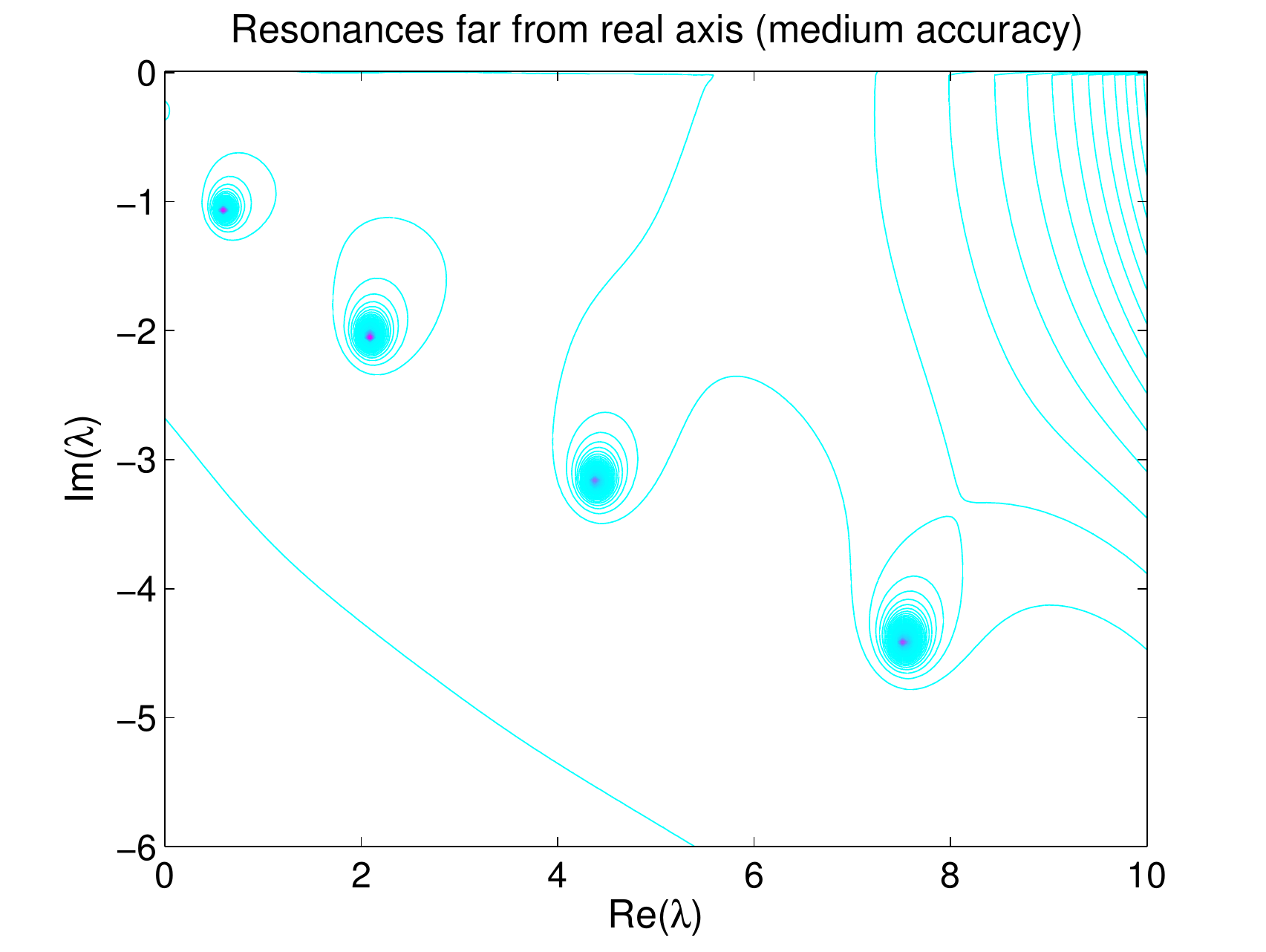}}\resizebox{250pt}{220pt}{\includegraphics*{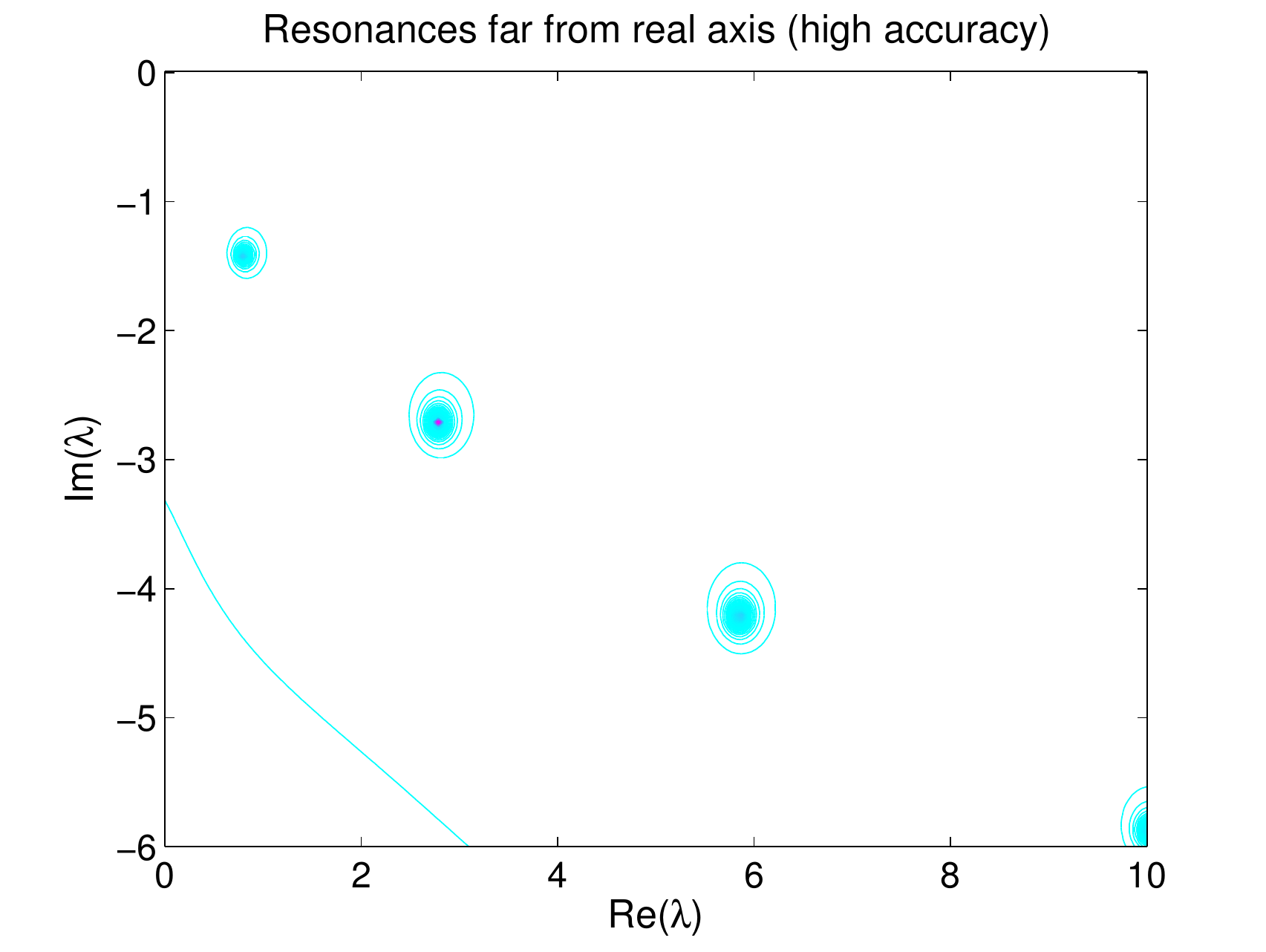}}
\caption{Resonances far from the real axis for $\epsilon=0.2$:
the medium and high accuracy results are only qualitatively similar\label{fig:mm4}}
\end{center}
\end{figure}

We also sought to calculate resonances further out in the complex plane.
These are rather more difficult to obtain --- see Fig. \ref{fig:mm4}, where
we compare the results for the medium and high levels of accuracy. The
agreement is rather poor. In fact the change in $R$ is more responsible
for this than the level of mesh refinement or the value of $\lambda_{\text{max}}$
for fixed $R$. This is rather different from the situation in \cite{MM},
where the resonances were calculated by complex scaling. There, the 
sensitive resonances depended more on the mesh than on the domain of
truncation.

\subsection{Example 4: Scattering by a Gaussian potential}
We consider the resonances of the Schr\"{o}dinger operator in $L_2 ( \Rbb[2] )$ given by
$-\Delta + q(x,y)$ where $q$ is a superposition of three Gaussians:
\[ q(x,y) = C\sum_{j=1}^{3}\exp(-\nu(x-x_j)^2-\nu(y-y_j)^2), \]
with $(x_1,y_1) = (0,-1)$, $(x_2,y_2) = (\sin(\pi/3),\cos(\pi/3))$,
$(x_3,y_3) = (-x_2,y_2)$, $C=40$ and $\nu=2$.  For large $C$ there is a high potential
barrier on the unit circle $x^2+y^2=1$ which gives rise to almost-trapped modes. However 
because $q$ is rapidly decaying the PDE $-\Delta u + qu = \lambda u$ in fact has no 
non-trivial $L_2$ solutions for $\lambda \geq 0$. The almost-trapped modes are the real
parts of resonances close to the real axis, which one may attempt to compute by the 
approach in Example 3.

An inner domain is created with boundary on the circle $x^2+y^2=R^2$. The value of
$R$ should be chosen so that for $x^2+y^2>R^2$, the function $q$ is small enough
to be neglected, which allows the PDE to be approximated by the Helmholtz equation. 
This allows us to use the same expression for the outer Neumann-to-Dirichlet map as
for Example 3.

As in Example 3, resonances far from the real axis are very unstable and cannot be
calculated reliably. However we were able to calculate some resonances close to the
real axis with moderate accuracy. One case is listed in Table \ref{table:mm4}. 
\begin{table}[tb!]
\begin{center}\small
\begin{tabular}{ccc}\hline
 number of mesh refines & $R$ & resonance \\ \hline
 2 &  4.0  & $4.257-5\times 10^{-4}\mathrm{i}$   \\
 3 &  4.0  & $4.140-5\times 10^{-4}\mathrm{i}$   \\
3 &  5.0  & $4.185-10^{-3}\mathrm{i}$           \\
4 &  6.0  & $4.402-4\times 10^{-4}\mathrm{i}$           \\
\hline
\end{tabular}
\end{center}
\caption{Resonance calculations for Gaussian potential. ($\lambda_{\text{max}}=50$)\label{table:mm4}}
\end{table}

Lin \cite{Lin} attempts to find these resonances using the complex scaling method and
a specially developed variational PDE-solver with a particularly chosen
set of basis functions. Even then, many spurious results are generated, and Lin 
claims to be able to filter them out.

\section{Conclusions}
Although our original motivation for doing this work lay in the need to devise
a numerical method which would be easy to implement, the numerical results show
that the accuracy achieved by our approach is actually quite competitive with
conventional techniques. Moreover, although Neumann series converge rather slowly,
two features ensure that for these problems the run-times were not excessive:
firstly, the fact that the quantities in the eqn. (\ref{eq:Rklinner}) 
are calculated once at the outset, so the most expensive calculations need
not be repeated for each different $\lambda$ at which the Neumann-to-Dirichlet
maps are evaluated; secondly, the trick in eqn. (\ref{eq:imprexpr}) acts as
a convergence acceleration technique for the Neumann expansion and means that 
in many cases we are able to get away with very few eigenfunctions. 

There are further modifications which could easily be made to improve
efficiency. Foremost among these would be to replace contour-plotting of
the condition number $\kappa(R(\lambda)+T(\lambda))$ by a routine
based on contour integration for finding zeros of an analytic function. There
are good practical reasons, however, for not doing this. Resonances close to
the real axis can already be located quite accurately by finding the maxima,
along the real axis, of $\kappa(R(\lambda)+T(\lambda))$; this
is easy since it is a one-dimensional search. As already observed by
Abramov, Aslanyan and Davies \cite{Ab}, resonances further from the real
axis are typically very unstable; to spend a lot of time accurately locating
the resonances of an approximating problem is therefore pointless.

\section*{Acknowledgements}\addcontentsline{toc}{section}{Acknowledgements} M.L. was partially supported by the Engineering and Physical Sciences Research Council 
grant EP/D054621/1. M.M. was partially supported by the Engineering and Physical Sciences Research Council grant EP/C008324/1.

\end{document}